\documentclass[table,svgnames,english,oneside]{smfart}

\title{Formalization of Harder-Narasimhan theory}
\usepackage{letltxmacro}
\LetLtxMacro\oldttfamily\ttfamily
\DeclareRobustCommand{\ttfamily}{\oldttfamily\csname ttsize\endcsname}
\newcommand{\setttsize}[1]{\def\ttsize{#1}}%
\usepackage{relsize}
\usepackage{orcidlink,amsmath,mathtools,amssymb,enumitem,geometry,graphicx,float}
\usepackage[format=plain,justification=centering]{caption}
\usepackage{mathrsfs}
\newcommand{\Memph}[1]{\textbf{\sffamily #1}}
\DeclareSymbolFont{bbold}{U}{bbold}{m}{n}
\DeclareMathSymbol{\bbmu}{\mathord}{bbold}{"16}
\usepackage{xparse}
\ExplSyntaxOn
\NewDocumentCommand{\definealphabet}{mmmm}
 {
  \int_step_inline:nnn { `#3 } { `#4 }
   {
    \cs_new_protected:cpx { #1 \char_generate:nn { ##1 }{ 11 } }
     {
      \exp_not:N #2 { \char_generate:nn { ##1 } { 11 } }
     }
   }
 }
\ExplSyntaxOff
\definealphabet{bb}{\mathbb}{A}{Z}
\definealphabet{cal}{\mathcal}{A}{Z}
\definealphabet{frak}{\mathfrak}{A}{z}
\definealphabet{rm}{\mathrm}{A}{z}
\definealphabet{bf}{\mathbf}{A}{Z}
\definealphabet{scr}{\mathscr}{A}{Z}
\definealphabet{tt}{\mathtt}{A}{Z}
\let\opn\operatorname

\usepackage{hyperref}
\usepackage[nameinlink]{cleveref}
\Crefname{theorem}{Theorem}{Theorem}
\Crefname{conjecture}{Conjecture}{Conjectures}
\Crefname{lemma}{Lemma}{Lemmas}
\Crefname{definition}{Definition}{Definitions}
\Crefname{remark}{Remark}{Remarks}
\Crefname{proposition}{Proposition}{Propositions}
\Crefname{corollary}{Corollary}{Corollaries}
\Crefname{equation}{}{}
\Crefname{item}{}{}
\Crefname{example}{Example}{Examples}
\Crefname{proof}{Proof}{Proofs}
\Crefname{condition}{Condition}{Conditions}
\Crefname{question}{Question}{Questions}
\newtheorem{theorem}{Theorem}[section]

\AddToHook{env/example/begin}{\crefalias{theorem}{example}}
\newtheorem{lemma}[theorem]{Lemma}
\AddToHook{env/lemma/begin}{\crefalias{theorem}{lemma}}
\newtheorem{remark}[theorem]{Remark}
\AddToHook{env/remark/begin}{\crefalias{theorem}{remark}}
\newtheorem{proposition}[theorem]{Proposition}
\AddToHook{env/proposition/begin}{\crefalias{theorem}{proposition}}
\newtheorem{definition}[theorem]{Definition}
\AddToHook{env/definition/begin}{\crefalias{theorem}{definition}}

\AddToHook{env/conjecture/begin}{\crefalias{theorem}{conjecture}}

\AddToHook{env/corollary/begin}{\crefalias{theorem}{corollary}}

\AddToHook{env/question/begin}{\crefalias{theorem}{question}}

\AddToHook{env/assumption/begin}{\crefalias{theorem}{assumption}}
\AddToHook{env/theorem/before}{\addvspace{3pt}}
\AddToHook{env/theorem/after}{\addvspace{3pt}}
\AddToHook{env/example/before}{\addvspace{3pt}}
\AddToHook{env/example/after}{\addvspace{3pt}}
\AddToHook{env/definition/before}{\addvspace{3pt}}
\AddToHook{env/definition/after}{\addvspace{3pt}}
\AddToHook{env/proposition/before}{\addvspace{3pt}}
\AddToHook{env/proposition/after}{\addvspace{3pt}}
\AddToHook{env/lemma/before}{\addvspace{3pt}}
\AddToHook{env/lemma/after}{\addvspace{3pt}}
\AddToHook{env/remark/before}{\addvspace{3pt}}
\AddToHook{env/remark/after}{\addvspace{3pt}}
\AddToHook{env/corollary/before}{\addvspace{3pt}}
\AddToHook{env/corollary/after}{\addvspace{3pt}}

\numberwithin{equation}{section}
\numberwithin{figure}{section}
\newlist{propenum}{enumerate}{1}
\setlist[propenum]{label=(\arabic*), ref=\theproposition~(\arabic*)}
\crefalias{propenumi}{proposition}
\newlist{lemenum}{enumerate}{1}
\setlist[lemenum]{label=(\arabic*), ref=\thelemma~(\arabic*)}
\crefalias{lemenumi}{lemma}
\newlist{thmenum}{enumerate}{1}
\setlist[thmenum]{label=(\arabic*), ref=\thetheorem~(\arabic*)}
\crefalias{thmenumi}{theorem}
\Crefformat{enumi}{#2\textup{(#1)}#3}
\newcommand{\link}[1]{\textsuperscript{\href{#1}{\faIcon{external-link-alt}}}}
\usepackage[url=true,backend=biber,style=ext-alphabetic,hyperref=true,giveninits=true,maxbibnames=99]{biblatex}
\usepackage{fontspec,fontawesome5}
\usepackage[dvipsnames]{xcolor}
\usepackage[frozencache]{minted}
\newmintinline[lean]{lean4}{bgcolor=white}
\usepackage{tcolorbox}
\tcbuselibrary{listings, minted, skins}

\tcbset{listing engine=minted}
\newtcblisting{leancode}[1][]{listing only, minted language=lean4, minted style=tango,
  colback=white!87!Emerald, enhanced, frame hidden, minted options={
  fontsize=\footnotesize, tabsize=2, breaklines, autogobble,mathescape},top=0pt,bottom=0pt,left=0pt,right=0pt,#1}

\author{Yijun Yuan\orcidlink{0000-0001-6571-6980}}
\address{Institute for Theoretical Sciences, Westlake University, No. 600 Dunyu Road, Sandun town, Xihu district, Hangzhou, Zhejiang Province, 310030, China}
\email{941201yuan@gmail.com}
\urladdr{https://yijunyuan.github.io/}
\begin{document}
\setttsize{\relscale{0.8}}

\begin{abstract}
  The Harder-Narasimhan theory provides a canonical filtration of a vector bundle on a projective curve whose successive quotients are semistable with strictly decreasing slopes. In this article, we present a formalization of Harder-Narasimhan theory in the proof assistant Lean 4 with Mathlib. The formalization is based on a recent approach to Harder-Narasimhan theory by Chen and Jeannin, which reinterprets the theory in order-theoretic terms and avoids the classical dependence on algebraic geometry. As an application, we formalize the uniqueness of the coprimary filtration of a nontrivial finitely generated module over a Noetherian ring, as well as the existence of a Jordan-Hölder filtration for a semistable Harder-Narasimhan game.
\end{abstract}
\subjclass{Primary 14H60, 68V20; Secondary 06B23, 91A80}
\keywords{Harder-Narasimhan theory, semistability, formalization, Lean 4, order theory}
\maketitle
\tableofcontents
\section{Introduction}
\subsection{Harder-Narasimhan theory}
Harder-Narasimhan theory, developed in the 1970s by Harder and Narasimhan, is a classical tool in algebraic geometry for studying vector bundles over algebraic curves through canonical filtrations with semistable successive quotients. We recall the relevant notions as follows:
\begin{definition}
    Let $X$ be a regular projective curve over a field $k$.
    \begin{enumerate}
        \item For any non-zero vector bundle $\calE$, define the \textbf{slope} of $\calE$, denoted by $\mu(\calE)$, as $\mu(\calE)\coloneqq\deg(\calE)/\opn{rk}(\calE)$, where $\deg(\calE)$ is the \textbf{degree} of $\calE$ and $\opn{rk}(\calE)$ is the \textbf{rank} of $\calE$.
        \item A vector bundle $\calE$ is called \textbf{semistable} if for every non-zero proper sub-bundle $\calF\subseteq\calE$, we have $\mu(\calF)\leq\mu(\calE)$.
    \end{enumerate}
\end{definition}
The main result of the Harder-Narasimhan theory is the following theorem:
\begin{theorem}[Harder-Narasimhan, {\cite{harderCohomologyGroupsModuli1975}}]\label{thm:classicalHN}
    Let $X$ be a regular projective curve over a field $k$. Any non-zero vector bundle $\calE$ over $X$ has a unique filtration by vector sub-bundles
    $$0=\calE_0\subsetneq\calE_1\subsetneq \cdots\subsetneq \calE_n=\calE,$$
    such that all subquotient sheaves $\calE_i/\calE_{i-1}$ are semistable vector bundles over $X$ and the successive slopes are strictly decreasing, i.e. $\mu(\calE_1/\calE_0)>\mu(\calE_2/\calE_1)>\cdots>\mu(\calE_n/\calE_{n-1})$.
\end{theorem}

This filtration, now known as the \textbf{Harder-Narasimhan filtration}, has been generalized to numerous contexts of rather different natures. For example,
\begin{enumerate}
    \item In \cite{hilleStableRepresentationsQuivers2002}, Hille and de la Peña generalized the Harder-Narasimhan filtration to the category of representations of quivers;
    \item In \cite{Fargues+2010+1+39}, Fargues studied a Harder-Narasimhan filtration of finite flat group schemes over an unequal characteristic valuation ring;
    \item In \cite{844f2731-ee1f-3791-921d-a8db397a2ba7}, Kedlaya introduced an analogue of the Harder-Narasimhan filtration for locally free sheaves on an open $p$-adic annulus equipped with a Frobenius structure to prove a $p$-adic monodromy theorem;
    \item In \cite{CML_2018__10_2_3_0}, Cornut used buildings to develop a metric version of the Harder-Narasimhan formalism on modular lattices of finite length, derived from it a categorical Harder-Narasimhan formalism, and gave a sufficient condition for compatibility with tensor products;
    \item In \cite{randriambololonaHarderNarasimhanTheory2019}, Randriambololona developed a Harder-Narasimhan theory for linear codes;
    \item In \cite{cornutHarderNarasimhanFiltrations2019}, Cornut and Peche Irissarry defined and studied a Harder-Narasimhan filtration on Breuil-Kisin-Fargues modules and related objects relevant to $p$-adic Hodge theory;
    \item In \cite{chenArakelovGeometryAdelic2020}, Chen and Moriwaki established a Harder-Narasimhan theory for adelic vector bundles.
\end{enumerate}
Notice that in \cite{chenArakelovGeometryAdelic2020}, for any non-zero normed lattice $(E,\lVert\cdot\rVert)$, the slope is defined as
$$\widehat{\mu}(E,\lVert\cdot\rVert)\coloneqq \frac{\widehat{\deg}(E,\lVert\cdot\rVert)}{\opn{rk}_{\bfZ}(E)},$$
where $\opn{rk}_{\bfZ}(E)$ is the rank of $E$ as a $\bfZ$-module and $\widehat{\deg}(E,\lVert\cdot\rVert)$ is the \textbf{Arakelov degree} of $(E,\lVert\cdot\rVert)$. Since the Arakelov degree is not additive with respect to short exact sequences, the traditional approach to establishing Harder-Narasimhan theory does not apply in this setting.

In view of the wide applicability of Harder-Narasimhan theory, and of the existence of variants, such as the one in \cite{chenArakelovGeometryAdelic2020}, that cannot be fully captured by the classical proof strategy, it is desirable to develop a more general yet tractable formulation that subsumes these generalizations.

In a recent work of Chen and Jeannin (cf. \cite{chen2023hardernarasimhangames}), Harder-Narasimhan theory is reinterpreted within order theory: the slope is interpreted as the pay-off function of a two-player zero-sum game on a bounded poset, and the semistability condition corresponds, under certain assumptions, to a Nash equilibrium of the game. This new perspective not only recovers the classical Harder-Narasimhan theory, but also applies to the case of normed lattices in \cite{chenArakelovGeometryAdelic2020}. The main result of \cite{chen2023hardernarasimhangames} is the following theorem:
\begin{theorem}[cf. {\cite[Theorem 1.1]{chen2023hardernarasimhangames}}]
  Let $(\scrL,\leq)$ be a bounded lattice, $(S,\leq)$ be a complete totally ordered set, and $\mu\colon \{(x,y)\in\scrL\times\scrL \vert x<y\}\to S$ be a map. Assume the following conditions:
\begin{enumerate}
  \item The pay-off function $\mu$ is convex (cf. \Cref{def:convex}).
  \item The bounded lattice $(\scrL,\leq)$ satisfies the ascending chain condition (i.e. no infinite strictly increasing sequence).
  \item $\mu$ satisfies the $\mu_A$-descending chain condition (cf. \Cref{def:dcc}).
\end{enumerate}
Then there exists a unique increasing sequence $\bot=a_0<a_1<\cdots<a_n=\top$
such that
\begin{enumerate}
  \item For any $i\in\{1,\cdots,n\}$, the restriction of the pay-off function $\mu$ to $\{z\in\scrL\vert a_{i-1}\leq z\leq a_i\}$ is semistable (cf. \Cref{def:semistable}).
  \item One has inequalities $\mu_A(\bot,a_1)>\mu_A(a_1,a_2)>\cdots>\mu_A(a_{n-1},\top)$, 
  where 
  $$\mu_A(x,y)\coloneqq\inf_{\substack{a\in\scrL\\x\leq a < y}}\sup_{\substack{b\in\scrL\\a<b\leq y}}\mu(a,b).$$
\end{enumerate}
\end{theorem}

In this article, we present the formalization of the version of Harder-Narasimhan theory due to Chen and Jeannin in the proof assistant Lean (cf. \cite{mouraLean4TheoremProver2021}, version 4.33.0) with Mathlib (cf. \cite{mathlib2020}, commit \href{https://github.com/leanprover-community/mathlib4/tree/db584cd6d46c92f209a44c0f1c829460d327499d}{db584cd}). Among the various versions of Harder-Narasimhan theory, we choose to formalize this one for the following reasons:
\begin{enumerate}
    \item \cite{chen2023hardernarasimhangames} strikes a good balance between generality and tractability. As mentioned above, the result in \cite{chen2023hardernarasimhangames} can be specialized to Harder-Narasimhan theory in various contexts, including the one in \cite{chenArakelovGeometryAdelic2020}. Meanwhile, compared with other attempts to categorify Harder-Narasimhan theory (e.g. \cite{chenHarderNarasimhanCategories2010}, \cite{CML_2018__10_2_3_0}, \cite{liCategorificationHarderNarasimhan2024} and \cite{pottharst2020hardernarasimhantheory}), the order-theoretic approach is easier to formalize.
    \item The exposition in \cite{chen2023hardernarasimhangames} is self-contained and modular, which greatly facilitates the formalization process.
    \item The Harder-Narasimhan theory in \cite{chen2023hardernarasimhangames} does not depend on any background from algebraic geometry, a subject that is not yet fully formalized in Mathlib.
\end{enumerate}
We refer the reader to \Cref{sec:mainthm} for the main theorem and its formalization.

\subsection{Related works}
Although this appears to be the first formalization of Harder-Narasimhan theory, it relies on Mathlib's order-theoretic and commutative-algebra foundations. We make the following incomplete list to record the contributions most relevant to this project.
\begin{enumerate}
  \item Leonardo de Moura formalized\link{https://github.com/leanprover-community/mathlib4/blob/db584cd6d46c92f209a44c0f1c829460d327499d/Mathlib/Order/Defs/PartialOrder.lean} the foundational order structures, including partial and linear orders.
  \item Johannes Hölzl formalized\link{https://github.com/leanprover-community/mathlib4/blob/db584cd6d46c92f209a44c0f1c829460d327499d/Mathlib/Order/CompleteLattice/Defs.lean}\link{https://github.com/leanprover-community/mathlib4/blob/db584cd6d46c92f209a44c0f1c829460d327499d/Mathlib/Order/BoundedOrder/Basic.lean} complete and bounded lattice foundations; Bhavik Mehta developed\link{https://github.com/leanprover-community/mathlib4/blob/db584cd6d46c92f209a44c0f1c829460d327499d/Mathlib/Order/Extension/Linear.lean} linear extensions of partial orders; Alena Gusakov and Yaël Dillies, together with Bhavik Mehta, developed\link{https://github.com/leanprover-community/mathlib4/blob/db584cd6d46c92f209a44c0f1c829460d327499d/Mathlib/Combinatorics/Colex.lean} colex orders; Violeta Hernández Palacios formalized\link{https://github.com/leanprover-community/mathlib4/blob/db584cd6d46c92f209a44c0f1c829460d327499d/Mathlib/Order/Completion.lean} the Dedekind--MacNeille completion; Jujian Zhang and Fangming Li formalized\link{https://github.com/leanprover-community/mathlib4/blob/db584cd6d46c92f209a44c0f1c829460d327499d/Mathlib/Order/RelSeries.lean} relation series; and Aaron Anderson and Yaël Dillies formalized\link{https://github.com/leanprover-community/mathlib4/blob/db584cd6d46c92f209a44c0f1c829460d327499d/Mathlib/Order/ModularLattice.lean} modular lattices.
  \item Chris Hughes formalized\link{https://github.com/leanprover-community/mathlib4/blob/db584cd6d46c92f209a44c0f1c829460d327499d/Mathlib/Order/JordanHolder.lean} the Jordan--Hölder theorem for abstract Jordan--Hölder lattices, including modular lattices, proving that any two composition series with the same endpoints are equivalent and hence have the same length.
  \item The commutative-algebra foundations of the coprimary-filtration application were developed by Jinzhao Pan (finiteness of the set of associated primes\link{https://github.com/leanprover-community/mathlib4/blob/db584cd6d46c92f209a44c0f1c829460d327499d/Mathlib/RingTheory/Ideal/AssociatedPrime/Finiteness.lean}), Nailin Guan (localization of associated primes\link{https://github.com/leanprover-community/mathlib4/blob/db584cd6d46c92f209a44c0f1c829460d327499d/Mathlib/RingTheory/Ideal/AssociatedPrime/Localization.lean}), Andrew Yang (localization of submodules\link{https://github.com/leanprover-community/mathlib4/blob/db584cd6d46c92f209a44c0f1c829460d327499d/Mathlib/Algebra/Module/LocalizedModule/Submodule.lean}), and Noah Walker (prime spectra\link{https://github.com/leanprover-community/mathlib4/blob/db584cd6d46c92f209a44c0f1c829460d327499d/Mathlib/RingTheory/Spectrum/Prime/Basic.lean}).
  \end{enumerate}

\begin{remark}
    When we started to formalize \cite{chen2023hardernarasimhangames}, the Dedekind-MacNeille completion was not yet a stable part of Mathlib, but only a draft pull request\link{https://github.com/leanprover-community/mathlib4/pull/26966}. We therefore had to implement it in our project\textsuperscript{\href{https://github.com/YijunYuan/HarderNarasimhan/blob/f0fd1d71851fe1d6aad52fa89446d02e0c26fdbd/HarderNarasimhan/OrderTheory/DedekindMacNeilleCompletion.lean}{\faIcon{external-link-alt}}}. After the pull request was merged, we removed our implementation and used the Mathlib version instead for consistency.
\end{remark}

Order theory, lattice theory, and commutative algebra have also been formalized in other proof assistants and their mathematical libraries. The following non-exhaustive list includes both foundational developments and more specialized results:
\begin{enumerate}
  \item Xavier Allamigeon, Cyril Cohen, Kazuhiko Sakaguchi, and Pierre-Yves Strub developed\link{https://coq-workshop.gitlab.io/2020/abstracts/Coq2020_03-02-ordered.pdf}\link{https://github.com/math-comp/math-comp/blob/98a2350f24de7fae13b5d54f592c8f835271cc69/order/order.v} a hierarchy of ordered structures in the Mathematical Components library for Coq (now Rocq), including partial and total orders, semilattices, lattices, and distributive lattices.
  \item John Harrison formalized\link{https://github.com/jrh13/hol-light/blob/cba9198db76e9dfb89cbd653df9412d01f65b22a/Library/wo.ml} foundational order theory in HOL Light, including partial and total orders, well-orders, and well-quasi-orders.
  \item John Harrison and Daniel J. Bernstein developed\link{https://github.com/jrh13/hol-light/blob/cba9198db76e9dfb89cbd653df9412d01f65b22a/Library/ringtheory.ml} a commutative-ring library in HOL Light, including ideals, quotient rings, localization, prime and maximal ideals, Noetherian rings, and polynomial rings.
  \item In Isabelle/HOL, Tobias Nipkow, Markus Wenzel, and Lawrence C. Paulson formalized\link{https://github.com/isabelle-prover/mirror-isabelle/blob/647253d043abe86e277658b8afb963b0ab76e491/src/HOL/Orderings.thy} preorders, partial orders, and linear orders; Nipkow developed\link{https://github.com/isabelle-prover/mirror-isabelle/blob/647253d043abe86e277658b8afb963b0ab76e491/src/HOL/Lattices.thy} abstract lattice theory; and Nipkow, Paulson, Wenzel, Florian Haftmann, and Viorel Preoteasa developed\link{https://github.com/isabelle-prover/mirror-isabelle/blob/647253d043abe86e277658b8afb963b0ab76e491/src/HOL/Complete_Lattices.thy} complete and completely distributive lattices.
  \item In Isabelle/HOL's HOL-Algebra library, Clemens Ballarin, with contributions from Martin Baillon and Jesús Aransay, developed\link{https://github.com/isabelle-prover/mirror-isabelle/blob/647253d043abe86e277658b8afb963b0ab76e491/src/HOL/Algebra/Ring.thy}\link{https://github.com/isabelle-prover/mirror-isabelle/blob/647253d043abe86e277658b8afb963b0ab76e491/src/HOL/Algebra/Module.thy}\link{https://github.com/isabelle-prover/mirror-isabelle/blob/647253d043abe86e277658b8afb963b0ab76e491/src/HOL/Algebra/UnivPoly.thy} rings, modules, and univariate polynomials; Stephan Hohe formalized\link{https://github.com/isabelle-prover/mirror-isabelle/blob/647253d043abe86e277658b8afb963b0ab76e491/src/HOL/Algebra/Ideal.thy} ideals, including prime and maximal ideals; Hohe and Paulo Emílio de Vilhena developed\link{https://github.com/isabelle-prover/mirror-isabelle/blob/647253d043abe86e277658b8afb963b0ab76e491/src/HOL/Algebra/QuotRing.thy} quotient rings; and de Vilhena formalized\link{https://github.com/isabelle-prover/mirror-isabelle/blob/647253d043abe86e277658b8afb963b0ab76e491/src/HOL/Algebra/Ring_Divisibility.thy} Noetherian rings and their characterization by the ascending chain condition on ideals.
  \item In IsarMathLib for Isabelle/ZF, Sławomir Kołodyński developed\link{https://github.com/SKolodynski/IsarMathLib/blob/e0d7f7d8aa523709e3b139305ba6570df7d68fdd/IsarMathLib/Order_ZF.thy}\link{https://github.com/SKolodynski/IsarMathLib/blob/e0d7f7d8aa523709e3b139305ba6570df7d68fdd/IsarMathLib/Lattice_ZF.thy} basic order theory, semilattices, and lattices; Daniel de la Concepción formalized\link{https://github.com/SKolodynski/IsarMathLib/blob/e0d7f7d8aa523709e3b139305ba6570df7d68fdd/IsarMathLib/Ring_ZF_2.thy}\link{https://github.com/SKolodynski/IsarMathLib/blob/e0d7f7d8aa523709e3b139305ba6570df7d68fdd/IsarMathLib/Ring_Zariski_ZF.thy} ring ideals, quotient rings, prime spectra, and the Zariski topology.
  \item In the Mizar Mathematical Library, Wojciech A. Trybulec formalized\link{https://mizar.uwb.edu.pl/version/current/html/orders_1.html} partially ordered sets; Stanisław Żukowski formalized\link{https://mizar.uwb.edu.pl/version/current/html/lattices.html} basic lattice theory, including modular lattices; and Grzegorz Bancerek formalized\link{https://mizar.uwb.edu.pl/version/current/html/lattice3.html} complete lattices.
  \item Jonathan Backer, Piotr Rudnicki, and Christoph Schwarzweller developed\link{https://mizar.uwb.edu.pl/version/current/html/ideal_1.html} ring ideal theory in Mizar, including generated ideals, ideal sums and radicals, 
  and Noetherian rings.
  \item Yasushige Watase developed\link{https://mizar.uwb.edu.pl/version/current/html/ideal_2.html} the theory of primary and quasi-primary ideals in Mizar.
  \item Joseph Hing-Lun Chan formalized\link{https://github.com/HOL-Theorem-Prover/HOL/blob/38cdeb6cdf10ad419dd81f730ac3eda8cb5d804c/src/algebra/construction/ringScript.sml} commutative ring theory in HOL4, including ring homomorphisms, ideals, quotient rings, 
  integral domains, etc.
  \item Ayberk Tosun formalized\link{https://github.com/martinescardo/TypeTopology/blob/19b7895d644def0dbbced98a097f7e4e3a46f633/source/Locales/DistributiveLattice/Definition.lagda}\link{https://github.com/martinescardo/TypeTopology/blob/19b7895d644def0dbbced98a097f7e4e3a46f633/source/Locales/DistributiveLattice/Ideal.lagda}\link{https://github.com/martinescardo/TypeTopology/blob/19b7895d644def0dbbced98a097f7e4e3a46f633/source/Locales/Frame.lagda} bounded distributive lattices, their ideals, and frames in Agda's TypeTopology library.
\end{enumerate}

\subsection{What we learned from this project}
\begin{enumerate}
  \item Formalization in Lean does not necessarily have to proceed strictly from the bottom up: even when the foundational machinery of a subject—such as the classical algebro-geometric setting of the Harder-Narasimhan theorem—is not yet fully available, it is still possible to formalize sophisticated results by identifying and working with suitable generalizations. By abstracting the essential ideas of Harder-Narasimhan theory beyond algebraic geometry, as in \cite{chen2023hardernarasimhangames}, we were able to capture its logical core in a setting where the necessary background, such as order theory, was already accessible. Once the algebraic geometry foundations in Mathlib become solid, the classical Harder-Narasimhan theorem can be recovered with little effort.
  
  This shows that creative reframing of mathematical concepts can open pathways for formalization long before all of the usual prerequisites are in place, and it encourages a flexible, problem-driven approach to building libraries in proof assistants.
  \item The process of formalization helps clarify arguments and reveal subtle gaps that often go unnoticed in conventional mathematical writing. While formalizing \cite{chen2023hardernarasimhangames}, we found that several necessary (but non-fatal) hypotheses were missing in some lemmas, such as\footnote{These minor gaps have been fixed in the latest version of \cite{chen2023hardernarasimhangames}.}:
  \begin{enumerate}
    \item In \cite[Proposition 4.8]{chen2023hardernarasimhangames}, the totally ordered $\bfR$-vector space $V$ needs to be non-trivial (which guarantees that the top element of \lean{DedekindCut V} is not a principal cut);
    \item In \cite[Theorem 4.21]{chen2023hardernarasimhangames}, to prove that condition (e) implies the others, one needs to assume the hypothesis of \cite[Proposition 4.20]{chen2023hardernarasimhangames};
    \item In \cite[Remark 4.26]{chen2023hardernarasimhangames}, one needs to add the condition that $\scrL$ is a modular lattice.
  \end{enumerate}
\end{enumerate}

\subsection{Future directions}
\begin{enumerate}
  \item Once the algebro-geometric framework within Mathlib is firmly established, it is expected to be straightforward to formalize the classical Harder-Narasimhan theorem (cf. \Cref{thm:classicalHN}) as a corollary of the main theorem in this project (cf. \Cref{thm:hnthm}).
  \item This project has also prompted us to reconsider, in the reverse direction, how to further develop the order-theoretic Harder-Narasimhan theory. For instance:
  \begin{itemize}
    \item Does this version of the Harder-Narasimhan filtration admit functoriality?
    \item Under what conditions is this filtration compatible with tensor products?
  \end{itemize}
    Answering these questions would deepen our understanding of the connections between the order-theoretic formulation of Harder-Narasimhan theory and the other formulations mentioned above.
\end{enumerate}

\subsection*{Acknowledgements}
The author would like to express his gratitude to Huayi Chen for his useful explanations about \cite{chen2023hardernarasimhangames}, to Jiedong Jiang for his helpful discussions about Lean, Mathlib and the related mathematical writing style, to Dingxin Zhang for his suggestion on typesetting, and to Junjie Bai, Zebei Li and Yicheng Tao for their help with Lean coding. Special thanks go to Shanwen Wang, who introduced the author to the field of formalized mathematics. Finally, the author would like to thank the anonymous referees for their valuable comments and suggestions, which have greatly improved the quality of this article.

The research is partially supported by the National Key R\&D Program of China (Grant No. 2024YFA1014000).

\subsection*{Code availability}The code of this project is available at \url{https://github.com/YijunYuan/HarderNarasimhan/tree/c6eb8349df607ee99fa2cab62aee7d9ec367418a}. Every code block in this article is clickable and hyperlinked to the corresponding code in the repository and Mathlib.

\section{Structure of the project}
The current library is organized around one infrastructure file and four functional blocks.
The umbrella module \lean{HarderNarasimhan.lean} imports the whole development.
\begin{enumerate}
    \item \lean{StrictIntvl.lean} defines strict intervals, their inclusion order, the total
    interval, and the type of points of an interval. The latter inherits the bounded-order,
    lattice, modularity, and well-foundedness instances needed by restrictions.
    \item \Memph{PayoffFunction/} contains the bundled pay-off function and its extremal
    operations, restriction and duality, convexity, semistability and breakpoints, slope-like
    functions, game values, and Nash equilibrium.
    \item \Memph{Filtration/} defines Harder-Narasimhan filtrations, constructs the canonical
    filtration by iterating greatest breakpoints, and proves uniqueness. It also provides the
    corresponding \lean{RelSeries} interface.
    \item \Memph{JordanHolder/} defines Jordan-H\"{o}lder filtrations, proves their existence and the stability of the restricted games on successive intervals, and proves equality of their lengths over a modular lattice.
    \item \Memph{Coprimary/} specializes the abstract theory to submodules of a finitely
    generated module over a Noetherian ring. It contains the required result on associated
    primes, the coprimary pay-off function and its semistability criterion, and the existence
    and uniqueness of coprimary filtrations.
\end{enumerate}

\begin{remark}
  We also provide a \href{https://github.com/YijunYuan/HarderNarasimhan/tree/literal}{\texttt{literal}} branch, which aims to follow the exposition and numbering of \cite{chen2023hardernarasimhangames} as closely as possible.
  This branch separates definitions, statements, and proofs into \lean{Defs.lean}, \lean{Results.lean}, and \lean{Impl.lean}, respectively, for each module, to facilitate reading the formalization alongside the paper.
\end{remark}

\section{The pay-off function and related conditions}
To introduce the central concept of Harder-Narasimhan theory, namely the semistability of the pay-off function, we start by recalling the definition of the pay-off function:
\begin{definition}\label{def:pay-off}
    Let $(\scrL,\leq)$ be a bounded poset with least element $\bot$ and greatest element $\top$, satisfying $\bot\neq \top$. Let $(S,\leq)$ be a complete lattice. 
	\begin{enumerate}
		\item A \textbf{strict interval} of $\scrL$ is a subset of $\scrL$ of the form $\scrL_{[x,y]}\coloneqq\{z\in\scrL|x\leq z \leq y\}$ for some $x,y\in\scrL$ with $x<y$. If no confusion arises, we simply call strict intervals \textbf{intervals}.
		\item The \textbf{pay-off function} $\mu$ of the \textbf{Harder-Narasimhan game} $(\scrL,S,\mu)$ is a map $$\mu\colon\{(x,y)\in \scrL^2:x<y\}\to S.$$
		In other words, the pay-off function assigns a value in $S$ to each strict interval of $\scrL$.
	\end{enumerate}
\end{definition}
Unless otherwise specified, $\scrL$, $S$, and $\mu$ are as in \Cref{def:pay-off} throughout the project.
\begin{remark}
    We do not distinguish between the properties of the Harder-Narasimhan game and the properties of its pay-off function. For example, when we say the game $(\scrL,S,\mu)$ is semistable, we mean the pay-off function $\mu$ is semistable.
\end{remark}
\begin{remark}
	At various places in the project, additional assumptions on $\scrL$ and $S$ are imposed:
\begin{enumerate}
    \item In the convexity, semistability, filtration, and Jordan-H\"older developments, we assume that $\scrL$ is a lattice and that $S$ is a complete lattice. Under these circumstances, the instance \lean{[PartialOrder ℒ]} is replaced by \lean{[Lattice ℒ]}, and the instance \lean{[CompleteLattice S]} is added.
    \item For some results (e.g. the uniqueness of the Harder-Narasimhan filtration), we need to assume that $(\scrL,\leq)$ satisfies the ascending chain condition (i.e. no infinite strictly increasing sequence). Under these circumstances, the instance \lean{[WellFoundedGT ℒ]} is added.
    \item For some results (e.g. the uniqueness of the Harder-Narasimhan filtration), we need to assume that $(S,\leq)$ is totally ordered. Under these circumstances, the instance \lean{[CompleteLinearOrder S]} is used instead of \lean{[CompleteLattice S]}.
\end{enumerate}
\end{remark}


\subsection{The pay-off function and its restriction \& dual}
Harder-Narasimhan theory studies not only the pay-off function of a fixed game, but also its restrictions to intervals. In this subsection, we present our design for the pay-off function and its restriction.

The domain of a pay-off function is represented by a dedicated structure:
\begin{leancode}[hyperurl={https://github.com/YijunYuan/HarderNarasimhan/blob/c6eb8349df607ee99fa2cab62aee7d9ec367418a/HarderNarasimhan/StrictIntvl.lean\#L38-L46}]
structure StrictIntvl (ℒ : Type*) [LT ℒ] where
  left : ℒ -- The left endpoint.
  right : ℒ -- The right endpoint.
  lt : left < right -- The endpoints are in strict order.
\end{leancode}
Such a strict interval packages its endpoints together with the proof that they are strictly
ordered. Strict intervals are ordered by inclusion, and when $\scrL$ is a nontrivial bounded
partial order, the total interval $(\bot,\top)$ is the top element of
\lean{StrictIntvl ℒ}. This is implemented as an \lean{OrderTop} instance of \lean{StrictIntvl ℒ}. 

In addition, we provide several helper lemmas in \lean{StrictIntvl.lean} to facilitate the manipulation of strict intervals. Among other things, we provide an \lean{instance} of \lean{Membership ℒ (StrictIntvl ℒ)}, which allows us to write \lean{x ∈ I} for a point \lean{x : ℒ} and a strict interval \lean{I : StrictIntvl ℒ}. In particular, rather than defining a second named
interval type, a \lean{CoeSort} instance makes the interval itself coerce to the subtype of
its points:
\begin{leancode}[hyperurl={https://github.com/YijunYuan/HarderNarasimhan/blob/c6eb8349df607ee99fa2cab62aee7d9ec367418a/HarderNarasimhan/StrictIntvl.lean\#L118-L120}]
instance [LT ℒ] [LE ℒ] : CoeSort (StrictIntvl ℒ) (Type _) :=
  ⟨fun I ↦ {x // x ∈ I}⟩
\end{leancode}
Thus \lean{↥I} is a nontrivial bounded order with endpoints \lean{I.left} and \lean{I.right}; lattice
and well-foundedness instances are inherited from the ambient order.  A strict interval of
\lean{↥I} is canonically coerced back to a strict interval of $\scrL$:
\begin{leancode}[hyperurl={https://github.com/YijunYuan/HarderNarasimhan/blob/c6eb8349df607ee99fa2cab62aee7d9ec367418a/HarderNarasimhan/StrictIntvl.lean\#L145-L151}]
def ofSub (J : StrictIntvl ↥I) : StrictIntvl ℒ :=
  ⟨J.left, J.right, Subtype.coe_lt_coe.2 J.lt⟩

instance : CoeOut (StrictIntvl ↥I) (StrictIntvl ℒ) := ⟨ofSub⟩
\end{leancode}

The pay-off function itself is bundled as follows:
\begin{leancode}[hyperurl={https://github.com/YijunYuan/HarderNarasimhan/blob/c6eb8349df607ee99fa2cab62aee7d9ec367418a/HarderNarasimhan/PayoffFunction/Defs.lean\#L40-L46}]
structure PayoffFunction (ℒ : Type*) [LT ℒ] (S : Type*) where
  toFun : StrictIntvl ℒ → S
\end{leancode}
    The \lean{FunLike} instance lets us apply $\mu$ directly to a strict interval.

\begin{remark}
    It is technically possible to define the pay-off function as a function of type \lean{ℒ → ℒ → S}. We instead use \lean{StrictIntvl ℒ} for the following reasons:
    \begin{enumerate}
        \item This definition is closer to the one in \cite{chen2023hardernarasimhangames}, which makes it easier to compare the formalization with the original text.
        \item It ensures that the pay-off function is applied only when $x<y$, while giving the endpoints and interval operations stable field names rather than repeatedly projecting from a subtype of pairs. For the same reason, we do not use \lean{Set.Icc} to represent the domain of the pay-off function.
        \item Its coercion to a sort identifies an interval $I$ with the subtype of its points. Consequently, restriction to $I$ is again a pay-off function on a self-contained bounded order.
    \end{enumerate}
\end{remark}

\begin{definition}
    The \textbf{restriction} of the Harder-Narasimhan game $(\scrL,S,\mu)$ to an interval $\scrL_{[x,y]}$ is the Harder-Narasimhan game $(\scrL_{[x,y]},S,\mu_{[x,y]})$, where $\mu_{[x,y]}$ is the restriction of $\mu$ to the domain $\scrL_{[x,y]}$.
\end{definition}
The restriction of the pay-off function is implemented as follows:
\begin{leancode}[hyperurl={https://github.com/YijunYuan/HarderNarasimhan/blob/c6eb8349df607ee99fa2cab62aee7d9ec367418a/HarderNarasimhan/PayoffFunction/Restrict.lean\#L30-L33}]
def restrict (μ : PayoffFunction ℒ S) (I : StrictIntvl ℒ) : PayoffFunction ↥I S :=
  ⟨fun J ↦ μ ↑J⟩
\end{leancode}
For example, for any \lean{I : StrictIntvl ℒ}, the semistability (cf. \Cref{def:semistable}) of the restriction of the pay-off function $\mu$ to $I$ is written as \lean{(μ.restrict I).IsSemistable}.

In addition, the symmetry between the two players allows us to define the \textbf{dual} of a pay-off function (resp. of a game), which is formalized as follows:
\begin{leancode}[hyperurl={https://github.com/YijunYuan/HarderNarasimhan/blob/c6eb8349df607ee99fa2cab62aee7d9ec367418a/HarderNarasimhan/PayoffFunction/Defs.lean\#L64-L69}]
def dual (μ : PayoffFunction ℒ S) : PayoffFunction ℒᵒᵈ Sᵒᵈ :=
  ⟨fun p ↦ OrderDual.toDual <| μ ⟨p.right, p.left, p.lt⟩⟩
\end{leancode}

\subsection{The $\mu$-series functions}
Several functions are associated with the pay-off function $\mu$; they facilitate the comparison with the classical Harder-Narasimhan theory:
\begin{definition}
    Let $\scrL$, $S$, and $\mu$ be as in \Cref{def:pay-off}. For any $x,y\in \scrL$ with $x<y$, we define
    $$\mu_{\max}(x,y)\coloneqq \sup_{\substack{w\in\scrL\\x<w\leq y}}\mu(x,w),\ \mu_{\min}(x,y)\coloneqq \inf_{\substack{w\in\scrL\\x\leq w<y}}\mu(w,y),$$
    $$\mu_A(x,y)\coloneqq\inf_{\substack{a\in\scrL\\x\leq a < y}}\mu_{\max}(a,y),\ \mu_B(x,y)\coloneqq\sup_{\substack{b\in\scrL\\x<b\leq y}}\mu_{\min}(x,b),$$
    $$\mu_A^*\coloneqq \mu_A(\bot,\top),\ \mu_B^*\coloneqq \mu_B(\bot,\top).$$
\end{definition}
\begin{remark}
    \cite{chen2023hardernarasimhangames} introduced two more functions, $\mu_{\max}(y)$ and $\mu_{\min}(y)$, which are nothing but $\mu_{\max}(\bot,y)$ and $\mu_{\min}(\bot,y)$, respectively. For simplicity, we do not formalize them.
\end{remark}
In the bundled interface these constructions are themselves pay-off functions, so they may
be evaluated on any strict interval. Their definitions in \lean{PayoffFunction/Defs.lean}
are:

\begin{leancode}[hyperurl={https://github.com/YijunYuan/HarderNarasimhan/blob/c6eb8349df607ee99fa2cab62aee7d9ec367418a/HarderNarasimhan/PayoffFunction/Defs.lean\#L81-L106}]
def max : PayoffFunction ℒ S :=
  ⟨fun I ↦ ⨆ (u : ℒ) (hu : u ∈ Set.Ioc I.left I.right), μ ⟨I.left, u, hu.1⟩⟩

def min : PayoffFunction ℒ S :=
  ⟨fun I ↦ ⨅ (u : ℒ) (hu : u ∈ Set.Ico I.left I.right), μ ⟨u, I.right, hu.2⟩⟩

def A : PayoffFunction ℒ S :=
  ⟨fun I ↦ ⨅ (a : ℒ) (ha : a ∈ Set.Ico I.left I.right), μ.max ⟨a, I.right, ha.2⟩⟩

def B : PayoffFunction ℒ S :=
  ⟨fun I ↦ ⨆ (b : ℒ) (hb : b ∈ Set.Ioc I.left I.right), μ.min ⟨I.left, b, hb.1⟩⟩
\end{leancode}

Note that all these definitions live in the namespace \lean{PayoffFunction}, so the dot notation allows us to write \lean{μ.max}, \lean{μ.min}, \lean{μ.A}, and \lean{μ.B} for the corresponding pay-off function \lean{μ}.
\begin{remark}
	In \cite{chen2023hardernarasimhangames}, the authors also define $\mu_A^*$ and $\mu_B^*$, which are the values of $\mu_A$ and $\mu_B$ on the total interval $\scrL=(\bot,\top)$. In the formalization, we do not define them separately, but instead use \lean{μ.A ⊤} and \lean{μ.B ⊤}, which are already simple enough, to represent them. We note that $\mu_A^*$ and $\mu_B^*$ are dual to each other in the following sense:
\begin{leancode}[hyperurl={https://github.com/YijunYuan/HarderNarasimhan/blob/c6eb8349df607ee99fa2cab62aee7d9ec367418a/HarderNarasimhan/PayoffFunction/GameValue.lean\#L169-L185}]
theorem A_top_dual : (μ.dual.A ⊤).ofDual = μ.B ⊤

theorem B_top_dual : (μ.dual.B ⊤).ofDual = μ.A ⊤
\end{leancode}
\end{remark}

\subsection{Conditions on the pay-off function}
\cite{chen2023hardernarasimhangames} introduced many conditions on the pay-off function $\mu$. In this subsection, we present those that are frequently used in the project.
\subsubsection{Convexity}
\begin{definition}\label{def:convex}
    Let $\scrL$ be a non-trivial bounded lattice. The pay-off function $\mu$ is
    \begin{enumerate}
        \item \textbf{convex} if the inequality
    $$\mu(x\wedge y,x)\leq \mu(y,x\vee y)$$
    holds for any $x,y\in \scrL$ with $x\not\leq y$.
        \item \textbf{affine} if the equality
    $$\mu(x\wedge y,x)= \mu(y,x\vee y)$$
    holds for any $x,y\in \scrL$ with $x\not\leq y$.
    \end{enumerate}
\end{definition}
\begin{leancode}[hyperurl={https://github.com/YijunYuan/HarderNarasimhan/blob/c6eb8349df607ee99fa2cab62aee7d9ec367418a/HarderNarasimhan/PayoffFunction/Convex.lean\#L52-L57}]
class IsConvex (μ : PayoffFunction ℒ S) : Prop where
  le : ∀ x y : ℒ, (h : ¬ x ≤ y) →
    μ ⟨x ⊓ y, x, inf_lt_left.2 h⟩ ≤ μ ⟨y, x ⊔ y, right_lt_sup.2 h⟩
\end{leancode}
\begin{leancode}[hyperurl={https://github.com/YijunYuan/HarderNarasimhan/blob/c6eb8349df607ee99fa2cab62aee7d9ec367418a/HarderNarasimhan/PayoffFunction/Convex.lean\#L94-L98}]
class IsAffine (μ : PayoffFunction ℒ S) : Prop where
  eq : ∀ x y : ℒ, (h : ¬ x ≤ y) →
    μ ⟨x ⊓ y, x, inf_lt_left.2 h⟩ = μ ⟨y, x ⊔ y, right_lt_sup.2 h⟩
\end{leancode}
In practice, one uses not only the global convexity of the pay-off function, but also its convexity on a strict interval. Although one can implement this as \lean{(μ.restrict I).IsConvex}, a flattened interface is provided for convenience:
\begin{leancode}[hyperurl={https://github.com/YijunYuan/HarderNarasimhan/blob/c6eb8349df607ee99fa2cab62aee7d9ec367418a/HarderNarasimhan/PayoffFunction/Convex.lean\#L59-L64}]
class IsConvexOn (μ : PayoffFunction ℒ S) (I : StrictIntvl ℒ) : Prop where
  le : ∀ x y : ℒ, x ∈ I → y ∈ I → (h : ¬ x ≤ y) →
    μ ⟨x ⊓ y, x, inf_lt_left.2 h⟩ ≤ μ ⟨y, x ⊔ y, right_lt_sup.2 h⟩
\end{leancode}
The following theorem shows that these two ways of expressing convexity on an interval are equivalent:
\begin{leancode}[hyperurl={https://github.com/YijunYuan/HarderNarasimhan/blob/c6eb8349df607ee99fa2cab62aee7d9ec367418a/HarderNarasimhan/PayoffFunction/Convex.lean\#L89-L92}]
theorem isConvexOn_iff_isConvex_restrict : μ.IsConvexOn I ↔ (μ.restrict I).IsConvex
\end{leancode}

Every affine pay-off function is convex, and affine pay-off functions remain affine after
restriction; both facts are registered as instances.
\subsubsection{Semistability}
\begin{definition}\label{def:semistable}
    Let $\scrL$ be a non-trivial bounded lattice. 
    \begin{enumerate}
        \item We say the pay-off function $\mu$ is \textbf{semistable} if for every $x>\bot$, $\mu_A(\bot,x)\not>\mu_A(\bot,\top)$.
        \item We say a semistable pay-off function $\mu$ is \textbf{stable} if for every $\bot<x<\top$, we have $\mu_A(\bot,x)\neq \mu_A(\bot,\top)$.
    \end{enumerate}
\end{definition}
These are properties of the bundled pay-off function:
\begin{leancode}[hyperurl={https://github.com/YijunYuan/HarderNarasimhan/blob/c6eb8349df607ee99fa2cab62aee7d9ec367418a/HarderNarasimhan/PayoffFunction/Semistable/Defs.lean\#L103-L113}]
class IsSemistable (μ : PayoffFunction ℒ S) : Prop where
  not_lt : ∀ x : ℒ, (hx : ⊥ < x) → ¬ μ.A ⊤ < μ.A ⟨⊥, x, hx⟩

class IsStable (μ : PayoffFunction ℒ S) : Prop extends μ.IsSemistable where
  ne : ∀ x : ℒ, (hx : ⊥ < x) → x < ⊤ → μ.A ⟨⊥, x, hx⟩ ≠ μ.A ⊤
\end{leancode}

\subsubsection{Chain conditions}
\begin{definition}[Descending chain conditions]\label{def:dcc}\leavevmode
    \begin{enumerate}
        \item We say that $\mu$ satisfies the \textbf{$\mu_A$-descending chain condition}
        if for any $a\in\scrL$, there does not exist any infinite descending chain
        $x_0>x_1>\cdots>x_n>x_{n+1}>\cdots>a$
        in $\scrL$ such that 
        $$\mu_A(a,x_0)<\mu_A(a,x_1)<\cdots<\mu_A(a,x_n)<\mu_A(a,x_{n+1})<\cdots.$$ 
		\begin{leancode}[hyperurl={https://github.com/YijunYuan/HarderNarasimhan/blob/c6eb8349df607ee99fa2cab62aee7d9ec367418a/HarderNarasimhan/PayoffFunction/Semistable/Defs.lean\#L51-L56}]
class ADCC (μ : PayoffFunction ℒ S) : Prop where
  dcc : ∀ a : ℒ, ∀ f : ℕ → ℒ, (h₁ : ∀ n : ℕ, f n > a) → StrictAnti f →
    ∃ N : ℕ, ¬ μ.A ⟨a, f N, h₁ N⟩ < μ.A ⟨a, f <| N + 1, h₁ <| N + 1⟩
		\end{leancode}
        \item The \textbf{strong descending chain condition} requires that for any descending chain $x_0>x_1>\cdots>x_n>x_{n+1}>\cdots$ in $\scrL$, there exists $N\in\bfN$ such that $\mu(\bot,x_N)\leq \mu(x_{N+1},x_N)$.
\begin{leancode}[hyperurl={https://github.com/YijunYuan/HarderNarasimhan/blob/c6eb8349df607ee99fa2cab62aee7d9ec367418a/HarderNarasimhan/PayoffFunction/GameValue.lean\#L63-L69}]
class StrongDCC (μ : PayoffFunction ℒ S) : Prop where
  exists_le : ∀ x : ℕ → ℒ, (saf : StrictAnti x) →
    ∃ N : ℕ, μ ⟨⊥, x N, lt_of_le_of_lt bot_le <| saf <| Nat.lt_add_one N⟩ ≤
      μ ⟨x (N+1), x N, saf <| Nat.lt_add_one N⟩
\end{leancode}
        \item The \textbf{eventually-$\top$ descending chain condition} is satisfied if for any descending chain $x_0>x_1>\cdots>x_n>x_{n+1}>\cdots$ in $\scrL$, there exists $N\in\bfN$ such that $\mu(x_{N+1},x_N)=+\infty$.
\begin{leancode}[hyperurl={https://github.com/YijunYuan/HarderNarasimhan/blob/c6eb8349df607ee99fa2cab62aee7d9ec367418a/HarderNarasimhan/JordanHolder/Defs.lean\#L63-L70}]
class EventuallyTopDCC (μ : PayoffFunction ℒ S) : Prop where
  exists_eq_top : ∀ x : ℕ → ℒ, (hx : StrictAnti x) →
    ∃ N : ℕ, μ ⟨x (N + 1), x N, hx (lt_add_one N)⟩ = ⊤
\end{leancode}
Note that this condition is stronger than the strong descending chain condition. The implication is formalized as an \lean{instance}.
    \end{enumerate}
\end{definition}

We also mention the following ascending chain condition:
\begin{definition}
	We say that $\mu$ satisfies the \textbf{weak ascending chain condition} if for any ascending chain $x_0<x_1<\cdots<x_n<x_{n+1}<\cdots$ in $\scrL$, there exists $N\in\bfN$ such that $\mu(x_N,x_{N+1})\leq \mu(x_N,\top)$.
\begin{leancode}[hyperurl={https://github.com/YijunYuan/HarderNarasimhan/blob/c6eb8349df607ee99fa2cab62aee7d9ec367418a/HarderNarasimhan/PayoffFunction/GameValue.lean\#L50-L56}]
class WeakACC (μ : PayoffFunction ℒ S) : Prop where
  exists_le : ∀ x : ℕ → ℒ, (smf : StrictMono x) →
    ∃ N : ℕ, μ ⟨x N, x (N+1), smf <| Nat.lt_add_one N⟩ ≤
      μ ⟨x N, ⊤, lt_of_lt_of_le (smf <| Nat.lt_add_one N) le_top⟩
\end{leancode}
\end{definition}
This condition holds trivially when $(\scrL,\leq)$ satisfies the ascending chain condition, i.e. when there is no infinite strictly increasing sequence in $\scrL$. The implication is formalized as an \lean{instance}.

\subsubsection{Slope-like conditions}\label{sec:slope-like}
The slope-like condition for pay-off functions captures the behavior of the slope function in the classical Harder-Narasimhan theory, i.e. the quotient of two additive functions satisfying certain properties.
\begin{definition}
    We say that the pay-off function $\mu$ is \textbf{slope-like} if for any elements $x<y<z$ in $\scrL$, the following statements hold:
    \begin{enumerate}
        \item Either $\mu(x,y)\leq \mu(x,z)$, or $\mu(y,z)<\mu(x,z)$.
        \item Either $\mu(x,y)< \mu(x,z)$, or $\mu(y,z)\leq\mu(x,z)$.
        \item Either $\mu(x,z)< \mu(x,y)$, or $\mu(x,z)\leq\mu(y,z)$.
        \item Either $\mu(x,z)\leq \mu(x,y)$, or $\mu(x,z)<\mu(y,z)$.
    \end{enumerate}
\end{definition}
\begin{leancode}[hyperurl={https://github.com/YijunYuan/HarderNarasimhan/blob/c6eb8349df607ee99fa2cab62aee7d9ec367418a/HarderNarasimhan/PayoffFunction/SlopeLike.lean\#L45-L57}]
class IsSlopeLike (μ : PayoffFunction ℒ S) : Prop where
  slopelike : ∀ (x y z : ℒ), (h : x < y ∧ y < z) →
    (μ ⟨x, y, h.1⟩ ≤ μ ⟨x, z, lt_trans h.1 h.2⟩ ∨ 
     μ ⟨y, z, h.2⟩ < μ ⟨x, z, lt_trans h.1 h.2⟩) ∧
    (μ ⟨x, y, h.1⟩ < μ ⟨x, z, lt_trans h.1 h.2⟩ ∨ 
     μ ⟨y, z, h.2⟩ ≤ μ ⟨x, z, lt_trans h.1 h.2⟩) ∧
    ... -- Assertions (3) and (4) omitted for brevity
\end{leancode}

The following theorem provides an important example of a slope-like pay-off function, which justifies the name of this condition.
\begin{theorem}[cf. {\cite[Proposition 4.8]{chen2023hardernarasimhangames}}]\label{thm:slope-like}
    Let $(V,\leq)$ be a non-trivial totally ordered vector space over $\bfR$, and $(S,\leq)$ be the smallest complete lattice containing $V$\footnote{This is called the Dedekind-MacNeille completion of $V$.}. Let 
    $$r\colon \{(x,y)\in\scrL\times\scrL | x < y\}\to \bfR_{\geq 0},\ d\colon \{(x,y)\in\scrL\times\scrL | x < y\}\to V$$
    be two maps that satisfy the following conditions:
    \begin{enumerate}
        \item For any elements $x<y<z$ in $\scrL$, one has
        $$d(x,z)=d(x,y)+d(y,z),\ r(x,z)=r(x,y)+r(y,z).$$
        \item For any elements $x<y$ in $\scrL$, if $r(x,y)=0$, then $d(x,y)>0$.
    \end{enumerate}
    Suppose that the pay-off function $\mu$ is given by
    $$\mu(x,y)\coloneqq \begin{cases}
        r(x,y)^{-1}d(x,y),& \text{ if }r(x,y)>0\\
        +\infty,& \text{ if }r(x,y)=0.
    \end{cases}$$
    Then $\mu$ is slope-like.
\end{theorem}
We first define the slope function for the two maps $r$ and $d$ without imposing any conditions on them:
\begin{leancode}[hyperurl={https://github.com/YijunYuan/HarderNarasimhan/blob/c6eb8349df607ee99fa2cab62aee7d9ec367418a/HarderNarasimhan/PayoffFunction/Slope.lean\#L50-L54}]
def slope (r : StrictIntvl ℒ → ℝ≥0) (d : StrictIntvl ℒ → V) :
    PayoffFunction ℒ (DedekindCut V) :=
  ⟨fun I ↦ if _ : 0 < r I then .principal ((r I)⁻¹ • d I) else ⊤⟩
\end{leancode}
Here \lean{DedekindCut V} is the Dedekind-MacNeille completion of $V$, and \lean{.principal} is the constructor of \lean{DedekindCut V} that takes an element of $V$ to the corresponding principal cut. In other words, \lean{.principal} is the canonical embedding of $V$ into its Dedekind-MacNeille completion.

We then formalize \Cref{thm:slope-like} as follows:
\begin{leancode}[hyperurl={https://github.com/YijunYuan/HarderNarasimhan/blob/c6eb8349df607ee99fa2cab62aee7d9ec367418a/HarderNarasimhan/PayoffFunction/Slope.lean\#L64-L130}]
theorem isSlopeLike_slope [Nontrivial V] (r : StrictIntvl ℒ → ℝ≥0) (d : StrictIntvl ℒ → V)
    (hd : ∀ (x y z : ℒ), (h₁ : x < y) → (h₂ : y < z) →
      d ⟨x, z, h₁.trans h₂⟩ = d ⟨x, y, h₁⟩ + d ⟨y, z, h₂⟩)
    (hr : ∀ (x y z : ℒ), (h₁ : x < y) → (h₂ : y < z) →
      r ⟨x, z, h₁.trans h₂⟩ = r ⟨x, y, h₁⟩ + r ⟨y, z, h₂⟩)
    (hpos : ∀ (x y : ℒ), (h : x < y) → r ⟨x, y, h⟩ = 0 → 0 < d ⟨x, y, h⟩) :
    (slope r d).IsSlopeLike
\end{leancode}

The slope-like condition admits an equivalent and more intuitive formulation:
\begin{lemma}[cf. {\cite[Proposition 4.6]{chen2023hardernarasimhangames}}]
  The pay-off function $\mu$ is slope-like if and only if it satisfies the following \textbf{seesaw property}: for any $x<y<z$ in $\scrL$, one and only one of the following three conditions is satisfied:
  \begin{enumerate}
    \item $\mu(x,y)<\mu(x,z)<\mu(y,z)$;
    \item $\mu(x,y)>\mu(x,z)>\mu(y,z)$;
    \item $\mu(x,y)=\mu(x,z)=\mu(y,z)$.
  \end{enumerate}
\begin{leancode}[hyperurl={https://github.com/YijunYuan/HarderNarasimhan/blob/c6eb8349df607ee99fa2cab62aee7d9ec367418a/HarderNarasimhan/PayoffFunction/SlopeLike.lean\#L65-L88}]
theorem isSlopeLike_iff_seesaw :
    μ.IsSlopeLike ↔ ∀ (x y z : ℒ), (h₁ : x < y) → (h₂ : y < z) →
      (μ ⟨x, y, h₁⟩ < μ ⟨x, z, h₁.trans h₂⟩ ∧ μ ⟨x, z, h₁.trans h₂⟩ < μ ⟨y, z, h₂⟩) ∨
      (μ ⟨x, z, h₁.trans h₂⟩ < μ ⟨x, y, h₁⟩ ∧ μ ⟨y, z, h₂⟩ < μ ⟨x, z, h₁.trans h₂⟩) ∨
      (μ ⟨x, y, h₁⟩ = μ ⟨x, z, h₁.trans h₂⟩ ∧ μ ⟨x, z, h₁.trans h₂⟩ = μ ⟨y, z, h₂⟩)
\end{leancode}
\end{lemma}

In practice, we find that the literal form of the seesaw property is less convenient to use than expected. For instance, suppose that we want to prove $\mu(x,y)>\mu(x,z)$ under the hypothesis $\mu(x,y)>\mu(y,z)$ for a slope-like $\mu$ with $x<y<z$. To invoke the seesaw property, we need to rule out the first case
\begin{leancode}[hyperurl={https://github.com/YijunYuan/HarderNarasimhan/blob/c6eb8349df607ee99fa2cab62aee7d9ec367418a/HarderNarasimhan/PayoffFunction/SlopeLike.lean\#L70}]
μ ⟨x, y, h₁⟩ < μ ⟨x, z, h₁.trans h₂⟩ ∧ μ ⟨x, z, h₁.trans h₂⟩ < μ ⟨y, z, h₂⟩
\end{leancode}
\noindent and the third case
\begin{leancode}[hyperurl={https://github.com/YijunYuan/HarderNarasimhan/blob/c6eb8349df607ee99fa2cab62aee7d9ec367418a/HarderNarasimhan/PayoffFunction/SlopeLike.lean\#L72}]
μ ⟨x, y, h₁⟩ = μ ⟨x, z, h₁.trans h₂⟩ ∧ μ ⟨x, z, h₁.trans h₂⟩ = μ ⟨y, z, h₂⟩
\end{leancode}
\noindent and then deconstruct the second case.

Our solution is to provide two \lean{iff} lemmas for each of the three cases, which allows us to prove all similar statements with at most two \lean{rw} tactic calls.
\begin{leancode}[hyperurl={https://github.com/YijunYuan/HarderNarasimhan/blob/c6eb8349df607ee99fa2cab62aee7d9ec367418a/HarderNarasimhan/PayoffFunction/SlopeLike.lean\#L103-L127}]
lemma IsSlopeLike.seesaw_total_lt_right_iff :
    μ ⟨x, z, h₁.trans h₂⟩ < μ ⟨y, z, h₂⟩ ↔ μ ⟨x, y, h₁⟩ < μ ⟨x, z, h₁.trans h₂⟩

lemma IsSlopeLike.seesaw_left_lt_right_iff :
    μ ⟨x, y, h₁⟩ < μ ⟨y, z, h₂⟩ ↔ μ ⟨x, y, h₁⟩ < μ ⟨x, z, h₁.trans h₂⟩

lemma IsSlopeLike.seesaw_right_lt_total_iff :
    μ ⟨y, z, h₂⟩ < μ ⟨x, z, h₁.trans h₂⟩ ↔ μ ⟨x, z, h₁.trans h₂⟩ < μ ⟨x, y, h₁⟩

... -- Other 3 lemmas omitted for brevity.
\end{leancode}

This refactoring follows a simple software engineering principle: the way a theorem is stated should reflect how it will be used. Providing small lemmas for rewriting avoids repeated case analysis in later proofs.

\subsubsection{Nash equilibrium}
In game theory, a Nash equilibrium is a situation in which no player can benefit from changing strategies while the other players keep theirs unchanged. The concept is used in \cite[Section 4.4]{chen2023hardernarasimhangames} as a new perspective for understanding the semistability condition of the Harder-Narasimhan game.
\begin{definition}
    The Harder-Narasimhan game $(\scrL,S,\mu)$ (i.e. the pay-off function $\mu$) has a \textbf{Nash equilibrium} if $\mu_A(\bot,\top)=\mu_B(\bot,\top)$.
\begin{leancode}[hyperurl={https://github.com/YijunYuan/HarderNarasimhan/blob/c6eb8349df607ee99fa2cab62aee7d9ec367418a/HarderNarasimhan/PayoffFunction/NashEquilibrium.lean\#L48-L53}]
class HasNashEquilibrium (μ : PayoffFunction ℒ S) : Prop where
  eq : μ.A ⊤ = μ.B ⊤
\end{leancode}
\end{definition}

One of the key observations in \cite{chen2023hardernarasimhangames} is that, under certain conditions, the semistability of the pay-off function $\mu$ is equivalent to the existence of a Nash equilibrium. This is the content of the following theorem:
\begin{theorem}[{cf. \cite[Theorem 4.21]{chen2023hardernarasimhangames}}]
    Let $(\scrL,\leq)$ be a non-trivial bounded lattice, and suppose that the order on $S$ is total. If the pay-off function $\mu$ is slope-like and satisfies both the weak ascending chain condition and the strong descending chain condition, then the following statements are equivalent:
    \begin{enumerate}
        \item $\mu_{\max}(\bot,\top)=\mu(\bot,\top)$;
        \item $\mu_{\min}(\bot,\top)=\mu(\bot,\top)$;
        \item $\mu_{\min}(\bot,\top)=\mu_{\max}(\bot,\top)$;
        \item $\mu$ has a Nash equilibrium.
    \end{enumerate}
    These conditions are implied by semistability. Conversely, the existence of a Nash equilibrium implies semistability if, for every $x\neq\bot$, the restriction $\mu_{[\bot,x]}$ satisfies the weak ascending chain condition and the weak slope-like condition at its top endpoint. For a slope-like pay-off function the latter condition is inherited by restriction.
\end{theorem}
The four-way equivalence is summarized in the theorem \lean{nashEquilibrium_tfae}:
\begin{leancode}[hyperurl={https://github.com/YijunYuan/HarderNarasimhan/blob/c6eb8349df607ee99fa2cab62aee7d9ec367418a/HarderNarasimhan/PayoffFunction/NashEquilibrium.lean\#L175-L183}]
theorem nashEquilibrium_tfae [μ.WeakACC] [μ.StrongDCC] :
    List.TFAE [μ.max ⊤ = μ ⊤, μ.min ⊤ = μ ⊤, μ.min ⊤ = μ.max ⊤, μ.HasNashEquilibrium]
\end{leancode}
The two directions relating equilibrium and semistability are exposed separately as
\begin{leancode}[hyperurl={https://github.com/YijunYuan/HarderNarasimhan/blob/c6eb8349df607ee99fa2cab62aee7d9ec367418a/HarderNarasimhan/PayoffFunction/NashEquilibrium.lean\#L201-L224}]
theorem IsSemistable.hasNashEquilibrium {S : Type*} [CompleteLinearOrder S]
    {μ : PayoffFunction ℒ S} (hμ : μ.IsSemistable) [μ.WeakACC] [μ.WeakSlopeLikeAtTop] :
    μ.HasNashEquilibrium

theorem isSemistable_of_hasNashEquilibrium {S : Type*} [CompleteLattice S]
    {μ : PayoffFunction ℒ S}
    (h₁ : ∀ x : ℒ, (hx : x ≠ ⊥) → 
      (μ.restrict ⟨⊥, x, bot_lt_iff_ne_bot.2 hx⟩).WeakACC)
    (h₂ : ∀ x : ℒ, (hx : x ≠ ⊥) →
      (μ.restrict ⟨⊥, x, bot_lt_iff_ne_bot.2 hx⟩).WeakSlopeLikeAtTop)
    (h : μ.HasNashEquilibrium) : μ.IsSemistable
\end{leancode}
This separation keeps the hypotheses of each implication explicit instead of packaging them into a large conjunction.

\section{Existence and uniqueness of Harder-Narasimhan filtration}
\subsection{Formalization of the filtration and the main theorem}\label{sec:mainthm}

We start with the definition of a Harder-Narasimhan filtration:
\begin{definition}\label{def:hnfil}
    Let $(\scrL,S,\mu)$ be a Harder-Narasimhan game, with $\scrL$ being a non-trivial bounded lattice. A \textbf{Harder-Narasimhan filtration} of $(\scrL,S,\mu)$ is an ascending finite filtration of $\scrL$
    $$\bot=a_0<a_1<\cdots<a_n=\top,$$
    such that
    \begin{enumerate}
        \item For any $j\in\{0,\cdots,n-1\}$, the restricted pay-off function $\mu_{[a_j,a_{j+1}]}$ is semistable;
        \item One has 
        $$\mu_A(a_0,a_1)\not\leq\mu_A(a_1,a_2)\not\leq \cdots\not\leq \mu_A(a_{n-1},a_n).$$
    \end{enumerate}
    We call $n$ the \textbf{length} of the Harder-Narasimhan filtration.
\end{definition}
\begin{leancode}[hyperurl={https://github.com/YijunYuan/HarderNarasimhan/blob/c6eb8349df607ee99fa2cab62aee7d9ec367418a/HarderNarasimhan/Filtration/Defs.lean\#L80-L105}]
structure HarderNarasimhanFiltration (μ : PayoffFunction ℒ S) where
  toFun : ℕ → ℒ -- The underlying chain; apply via the coercion
  length : ℕ -- The index at which the chain reaches `⊤`.
  monotone : Monotone toFun
  head_eq_bot : toFun 0 = ⊥
  length_eq_top : toFun length = ⊤
  strictMonoOn : StrictMonoOn toFun (Set.Iic length)
  piecewise_isSemistable : ∀ i, (hi : i < length) →
    (μ.restrict ⟨toFun i, toFun (i + 1), ...⟩).IsSemistable
  not_A_le_succ : ∀ i, (hi : i + 1 < length) →
    ¬ μ.A ⟨toFun i, toFun (i + 1), ...⟩ ≤ μ.A ⟨toFun (i + 1), toFun (i + 2), ...⟩
\end{leancode}
\begin{remark}
  An alternative way to express that a filtration has finite length in Lean, used in an earlier version of this project (\href{https://github.com/YijunYuan/HarderNarasimhan/tree/0f2089dd193ee9a9d7f02f233921560e642a4a69}{commit 0f2089d}), is to include the field
\begin{leancode}[hyperurl={https://github.com/YijunYuan/HarderNarasimhan/blob/0f2089dd193ee9a9d7f02f233921560e642a4a69/HarderNarasimhan/Filtration/Defs.lean\#L95}]
fin_len : ∃ n : ℕ, toFun n = ⊤
\end{leancode}
\noindent in place of an explicit \lean{length} field. The length is then recovered using \lean{Nat.find fin_len}. Although this formulation closely reflects the mathematical statement that ``the filtration has finite length'', it makes subsequent use of the structure more verbose.
\end{remark}
Now we are ready to present the main theorem:
\begin{theorem}[cf. {\cite[Definition 3.9, Theorem 3.10]{chen2023hardernarasimhangames}}]\label{thm:hnthm}
    Let $(\scrL,S,\mu)$ be a Harder-Narasimhan game as in \Cref{def:hnfil}. Then $(\scrL,S,\mu)$ admits a Harder-Narasimhan filtration if the following conditions hold:
    \begin{enumerate}
        \item $(\scrL,\leq)$ satisfies the ascending chain condition;
        \item The pay-off function $\mu$ is convex and satisfies the $\mu_A$-descending chain condition;
        \item (\textbf{$\mu$-admissible condition}) $S$ is totally ordered, or the infimum
        $$\mu_A(x,y)=\inf_{\substack{a\in \scrL\\x\leq a<y}}\mu_{\max}(a,y)$$
        is attained for any elements $x<y$ in $\scrL$.
    \end{enumerate}
    Moreover, if $S$ is totally ordered, then the Harder-Narasimhan filtration is unique.
\end{theorem}
The admissibility hypothesis is itself bundled, and we provide an \lean{instance} for the case of a totally ordered $S$:
\begin{leancode}[hyperurl={https://github.com/YijunYuan/HarderNarasimhan/blob/c6eb8349df607ee99fa2cab62aee7d9ec367418a/HarderNarasimhan/Filtration/Defs.lean\#L54-L72}]
class Admissible (μ : PayoffFunction ℒ S) : Prop where
  total_or_attained : Std.Total (· ≤ · : S → S → Prop) ∨ ∀ I : StrictIntvl ℒ, μ.IsAttained I

instance [CompleteLinearOrder S] (μ : PayoffFunction ℒ S) : μ.Admissible
\end{leancode}
Under \lean{[WellFoundedGT ℒ]}, \lean{[μ.ADCC]}, \lean{[μ.IsConvex]}, and
\lean{[μ.Admissible]}, the public construction \lean{μ.hnFiltration} is an inhabitant of
\lean{μ.HarderNarasimhanFiltration}. When $S$ is a complete linear order, the
\lean{Unique (μ.HarderNarasimhanFiltration)} instance proves that every such filtration is
this canonical one.

The project also exposes a \lean{RelSeries} presentation through the relation
\lean{μ.semistableRel}:
\begin{leancode}[hyperurl={https://github.com/YijunYuan/HarderNarasimhan/blob/c6eb8349df607ee99fa2cab62aee7d9ec367418a/HarderNarasimhan/Filtration/Defs.lean\#L165-L168}]
def semistableRel (μ : PayoffFunction ℒ S) : SetRel ℒ ℒ :=
  {(x, y) | ∃ h : x < y, (μ.restrict ⟨x, y, h⟩).IsSemistable}
\end{leancode}
Existence and uniqueness in this form are stated by
\lean{exists_relSeries_semistableRel} and
\begin{leancode}[hyperurl={https://github.com/YijunYuan/HarderNarasimhan/blob/c6eb8349df607ee99fa2cab62aee7d9ec367418a/HarderNarasimhan/Filtration/Unique.lean\#L231-L260}]
theorem existsUnique_relSeries_semistableRel (μ : PayoffFunction ℒ S)
    [μ.ADCC] [μ.IsConvex] :
    ∃! s : RelSeries μ.semistableRel,
      s.head = ⊥ ∧ s.last = ⊤ ∧
      ∀ i : ℕ, (hi : i + 1 < s.length) →
        ¬ μ.A ⟨s ↑i, s ↑(i + 1), ...⟩ ≤ μ.A ⟨s ↑(i + 1), s ↑(i + 2), ...⟩
\end{leancode}
\noindent respectively. 

The structure API is used for
the canonical construction, while the \lean{RelSeries} API provides interoperability with
Mathlib.

\subsection{Proof sketch of the main theorem}
Note that \Cref{def:hnfil} is not equivalent to \cite[Definition 3.9]{chen2023hardernarasimhangames}, where the authors constructed a canonical filtration $\calF_{\text{can}}$ of $\scrL$ (under the assumptions of \Cref{thm:hnthm}) that satisfies the conditions in \Cref{def:hnfil}, and called it the Harder-Narasimhan filtration. If $S$ is not totally ordered, then there may exist another filtration $\calF'$ that satisfies the conditions in \Cref{def:hnfil}, and $\calF'$ is not the Harder-Narasimhan filtration in the sense of \cite[Definition 3.9]{chen2023hardernarasimhangames}.

In this project, we adopt the more descriptive \Cref{def:hnfil}, as it offers better intuition and aligns more closely with the formulation of \Cref{thm:classicalHN}. Nevertheless, the existence of the Harder-Narasimhan filtration is proved by constructing the filtration $\calF_{\text{can}}$ as in \cite[Definition 3.9]{chen2023hardernarasimhangames}.

Since the proof of \Cref{thm:hnthm} in \cite{chen2023hardernarasimhangames} is clear and concise, we only briefly sketch the construction of $\calF_{\text{can}}$ here. The following result is crucial in the construction and the uniqueness proof:
\begin{theorem}[cf. {\cite[Proposition 3.4, Remark 3.5, Definition 3.6, Proposition 3.8]{chen2023hardernarasimhangames}}]\label{thm:hnprop}
    Let $(\scrL,S,\mu)$ be a Harder-Narasimhan game, with $\scrL$ being a non-trivial bounded lattice and $\mu$ being a convex pay-off function. Let $\opn{St}(\mu)$ be the set of elements $x$ in $\scrL\backslash\{\bot\}$ that satisfy the following conditions:
    \begin{enumerate}
        \item for any $y\in\scrL\backslash\{\bot\}$, $\mu_A(\bot,y)\not>\mu_A(\bot,x)$;
        \item for any $y\in\scrL\backslash\{\bot\}$, if $\mu_A(\bot,y)=\mu_A(\bot,x)$, then $y\leq x$.
    \end{enumerate}
    Then
    \begin{enumerate}
        \item If $(\scrL,\leq)$ satisfies the ascending chain condition and $\mu$ satisfies the $\mu_A$-descending chain condition, then $\opn{St}(\mu)$ is non-empty.
        \item $\mu$ is semistable if and only if $\top\in \opn{St}(\mu)$. In particular, when $(S,\leq)$ is totally ordered, $\mu$ is semistable if and only if $\opn{St}(\mu)=\{\top\}$.
        \item For any $x\in \opn{St}(\mu)\backslash\{\top\}$,
        \begin{enumerate}
            \item the restricted pay-off function $\mu_{[\bot,x]}$ is semistable;
            \item for any $y\in\scrL$ such that $y>x$, one has $\mu_A(\bot,x)\not\leq \mu_A(x,y)$.
        \end{enumerate}
        \item If conditions of \Cref{thm:hnthm} hold, then $\opn{St}(\mu)$ admits a greatest element.
    \end{enumerate}
\end{theorem}
With this theorem, one can construct the canonical filtration $\calF_{\text{can}}$ of $\scrL$ as follows:
\begin{equation}\label{eq:40661}
	\forall n\in \bfN,\ \calF_{\text{can}}(n) \coloneqq \begin{cases}
    \bot,& \text{ if }n=0;\\
    \max \opn{St}\left(\mu_{[\calF_{\text{can}}(n-1),\top]}\right),& \text{ if } n>0 \text{ and } \calF_{\text{can}}(n-1)\neq \top;\\
    \top,& \text{ if } n>0 \text{ and } \calF_{\text{can}}(n-1) = \top.
\end{cases}
\end{equation}

It follows immediately from \Cref{thm:hnprop} that $\calF_{\text{can}}$ is indeed a Harder-Narasimhan filtration of $(\scrL,S,\mu)$ in the sense of \Cref{def:hnfil}.

\subsubsection*{Implementation}
Note that the pay-off function used in the set $\opn{St}(\cdot)$ in \eqref{eq:40661} varies with $n$. Instead of using \lean{restrict} here, we flatten the definition of $\opn{St}(\cdot)$ on a strict interval \texttt{I}:
\begin{leancode}[hyperurl={https://github.com/YijunYuan/HarderNarasimhan/blob/c6eb8349df607ee99fa2cab62aee7d9ec367418a/HarderNarasimhan/PayoffFunction/Semistable/Defs.lean\#L66-L87}]
structure IsBreakpoint (x : ℒ) : Prop where -- A structure-shaped predicate for $\operatorname{St}$
  mem : x ∈ I
  ne_left : I.left ≠ x
  not_lt : ∀ y : ℒ, (hyI : y ∈ I) → (hy : I.left ≠ y) →
    ¬ μ.A ⟨I.left, x, ...⟩ < μ.A ⟨I.left, y, ...⟩
  le_of_eq : ∀ y : ℒ, (hyI : y ∈ I) → (hy : I.left ≠ y) →
    μ.A ⟨I.left, y, ...⟩ = μ.A ⟨I.left, x, ...⟩ → y ≤ x

def breakpoints : Set ℒ := {x | μ.IsBreakpoint I x} -- The set $\opn{St}(\mu|_I)$
\end{leancode}

Then the relative-to-$I$ version of \Cref{thm:hnprop} (4) is formalized as follows:
\begin{leancode}[hyperurl={https://github.com/YijunYuan/HarderNarasimhan/blob/c6eb8349df607ee99fa2cab62aee7d9ec367418a/HarderNarasimhan/PayoffFunction/Semistable/Breakpoints.lean\#L351-L362}]
lemma exists_isGreatest_breakpoints [hwf : WellFoundedGT ℒ] [μ.ADCC] (hμcvx : μ.IsConvexOn I)
    (h : (Std.Total (· ≤ · : S → S → Prop)) ∨
      ∀ z : ℒ, (hzI : z ∈ I) → (hz : I.left ≠ z) →
        μ.IsAttained ⟨I.left, z, lt_of_le_of_ne hzI.left hz⟩) :
    ∃ s : ℒ, IsGreatest (μ.breakpoints I) s
\end{leancode}
This allows us to construct the canonical filtration $\calF_{\text{can}}$ recursively as a raw function from $\bfN$ to $\scrL$:
\begin{leancode}[hyperurl={https://github.com/YijunYuan/HarderNarasimhan/blob/c6eb8349df607ee99fa2cab62aee7d9ec367418a/HarderNarasimhan/Filtration/Exists.lean\#L50-L63}]
def HNFil (k : ℕ) : ℒ :=
  match k with
  | 0 => ⊥
  | n + 1 =>
    let prev := HNFil n
    if htop : prev = ⊤ then ⊤
    else
      (exists_isGreatest_breakpoints (I := ⟨prev, ⊤, lt_top_iff_ne_top.2 htop⟩)
        ... ...).choose
\end{leancode}
With this \lean{HNFil} constructed, we establish several formalized lemmas to verify that it satisfies the conditions in \Cref{def:hnfil}, and assemble them into \lean{μ.hnFiltration}:
\begin{leancode}[hyperurl={https://github.com/YijunYuan/HarderNarasimhan/blob/c6eb8349df607ee99fa2cab62aee7d9ec367418a/HarderNarasimhan/Filtration/Exists.lean\#L144-L159}]
def hnFiltration : μ.HarderNarasimhanFiltration where
  toFun := HNFil μ
  length := HNlen μ
  monotone := HNFil_monotone μ
  head_eq_bot := rfl
  ... -- Omitted for brevity.
\end{leancode}

\section{Application 1: the coprimary filtration}
\subsection{Prime filtration and coprimary filtration}
The prime filtration of a finitely generated module over a Noetherian ring is a classical result in commutative algebra, which we now briefly review.
\begin{definition}
    Let $R$ be a commutative ring and let $M$ be a module over $R$. 
    \begin{enumerate}
        \item Denote by $\opn{Ass}(M)$ the set of \textbf{associated primes} of $M$, i.e. the primes that arise as the annihilator of some cyclic submodule of $M$.
        \item Call 
        $$\opn{Supp}(M)\coloneqq \{\frakp\in \opn{Spec}(R)| M_{\frakp}\neq 0\}$$
        the \textbf{support} of $M$.
        \item If $\opn{Ass}(M)$ consists of a single prime $\frakp$, we call $M$ a ($\frakp$-)\textbf{coprimary} module.
\begin{leancode}[hyperurl={https://github.com/YijunYuan/HarderNarasimhan/blob/c6eb8349df607ee99fa2cab62aee7d9ec367418a/HarderNarasimhan/Coprimary/Defs.lean\#L115-L118}]
class IsCoprimary : Prop where
  existsUnique_associatedPrime : ∃! p, p ∈ associatedPrimes R M
\end{leancode}
    \end{enumerate}
\end{definition}
\begin{theorem}[cf. {\cite[Ch IV, §1, n°4, Théorème 1]{bourbakiAlgebreCommutative2006a}}]
    Let $R$ be a commutative Noetherian ring, and let $M$ be a non-trivial finitely generated $R$-module.
    \begin{enumerate}
        \item There exists a filtration of $M$ by submodules
        $$0 = M_0 \subsetneq M_1 \subsetneq \cdots \subsetneq M_n = M$$
        such that one has $M_{i+1}/M_i\cong R/\frakp_i$ for some prime ideal $\frakp_i$ of $R$.
        \item One has the following inclusion:
        \begin{equation}\label{eq:pfinc}
            \opn{Ass}(M)\subseteq\{\frakp_0,\cdots,\frakp_{n-1}\}\subseteq \opn{Supp}(M).
        \end{equation}
    \end{enumerate}
\end{theorem}
As noted in \cite[Ch IV, §1, n°4]{bourbakiAlgebreCommutative2006a}, the prime filtration is generally not unique, and the inclusion 
\begin{equation}\label{eq:ass}
    \opn{Ass}(M)\subseteq \{\frakp_0,\cdots,\frakp_{n-1}\}
\end{equation}
in \eqref{eq:pfinc} is not necessarily an equality.

On the other hand, \cite[Ch IV, §2, Exercise (7)]{bourbakiAlgebreCommutative2006a} shows that \Cref{eq:ass} can be made into an equality if one relaxes the condition $M_{i+1}/M_i\cong R/\frakp_i$ to the condition that $M_{i+1}/M_i$ is coprimary. The resulting filtration is called a \textbf{coprimary filtration of $M$}.

The new result in \cite[Section 3.4]{chen2023hardernarasimhangames} shows that, if one interprets the coprimary filtration as the Harder-Narasimhan filtration of a suitably constructed Harder-Narasimhan game, then the coprimary filtration is unique once a linear order on $\opn{Spec}(R)$ extending the partial order of inclusion is fixed:
\begin{theorem}[cf. {\cite[Theorem 3.15, Remark 3.16]{chen2023hardernarasimhangames}}]\label{thm:cofil}
    Let $M\neq 0$ be a finitely generated module over a commutative Noetherian ring $R$. We equip $\opn{Spec}(R)$ with a linear order that extends\footnote{Such a \textbf{linear extension} always exists by the Szpilrajn extension theorem (cf. \cite{szpilrajnLextensionLordrePartiel1930}).} the order given by inclusion. Then there exists a unique coprimary filtration of $M$
        $$0 = M_0 \subsetneq M_1 \subsetneq \cdots \subsetneq M_n = M$$
        such that each quotient $M_{i+1}/M_i$ is a $\frakp_i$-coprimary module, and $\frakp_0>\cdots>\frakp_{n-1}$. Moreover, one has
        $$\opn{Ass}(M)=\{\frakp_0,\cdots,\frakp_{n-1}\}.$$
\end{theorem}
In parallel with the Harder-Narasimhan filtration, we formalize the coprimary filtration as a \lean{structure}:
\begin{leancode}[hyperurl={https://github.com/YijunYuan/HarderNarasimhan/blob/c6eb8349df607ee99fa2cab62aee7d9ec367418a/HarderNarasimhan/Coprimary/Defs.lean\#L124-L153}]
structure CoprimaryFiltration (R : Type*) [CommRing R] [IsNoetherianRing R]
    (M : Type*) [Nontrivial M] [AddCommGroup M] [Module R M] [Module.Finite R M] where
  toFun : ℕ → Submodule R M
  length : ℕ
  monotone : Monotone toFun
  head_eq_bot : toFun 0 = ⊥
  length_eq_top : toFun length = ⊤
  strictMonoOn : StrictMonoOn toFun (Set.Iic length)
  /-- Each successive subquotient `F (i + 1) ⧸ F i` is coprimary. -/
  piecewise_isCoprimary : ∀ i < length,
    IsCoprimary R (toFun (i + 1) ⧸ (toFun i).submoduleOf (toFun (i + 1)))
  /-- The associated primes of the successive subquotients strictly decrease along the
  chain, in the fixed linear extension of the prime spectrum. -/
  associatedPrime_succ_lt : ∀ i, i + 1 < length → ∀ p q : PrimeSpectrum R,
    p.asIdeal ∈ associatedPrimes R
      (toFun (i + 2) ⧸ (toFun (i + 1)).submoduleOf (toFun (i + 2))) →
    q.asIdeal ∈ associatedPrimes R
      (toFun (i + 1) ⧸ (toFun i).submoduleOf (toFun (i + 1))) →
    toLinearExtension p < toLinearExtension q
\end{leancode}
\noindent and construct the canonical coprimary filtration under the name \lean{Coprimary.coprimaryFiltration}. The uniqueness is expressed by the \lean{Unique (CoprimaryFiltration R M)} instance, and the final assertion of \Cref{thm:cofil} is formalized as
\begin{leancode}[hyperurl={https://github.com/YijunYuan/HarderNarasimhan/blob/c6eb8349df607ee99fa2cab62aee7d9ec367418a/HarderNarasimhan/Coprimary/Filtration.lean\#L143-L194}]
theorem associatedPrimes_eq_iUnion (F : CoprimaryFiltration R M) :
    associatedPrimes R M =
      ⋃ i < F.length, associatedPrimes R (F (i + 1) ⧸ (F i).submoduleOf (F (i + 1)))
\end{leancode}

\subsection{Realizing coprimary filtration as Harder-Narasimhan filtration}
As mentioned above, \Cref{thm:cofil} is proved by realizing coprimary filtrations as the Harder-Narasimhan
filtrations of a suitable bundled pay-off function and then applying \Cref{thm:hnthm}.

\subsubsection{The bounded lattice $\scrL_M$}
Let $\scrL_M$ be the lattice \lean{Submodule R M}.  Its order is inclusion, its join is the
sum of submodules, and its meet is their intersection.  All bounded-lattice and complete-lattice
instances are supplied by Mathlib.  Since $M$ is finite over the Noetherian ring $R$, this
lattice also satisfies the ascending chain condition.  The current implementation uses
\lean{Submodule R M} directly, rather than introducing a project-wide abbreviation.

\subsubsection{The complete lattice $S_R$}
The required complete lattice $S_R$ is built from the set of all finite subsets of the prime spectrum $\opn{Spec}(R)$, as follows.

To enrich the order structure on $\opn{Spec}(R)$, the raw set \lean{PrimeSpectrum R} is first replaced by \lean{LinearExtension} applied to
it, which fixes a linear order on $\opn{Spec}(R)$ that extends the partial order of inclusion. We then equip the set of all finite sets of primes with a total order $\leq$ that extends the partial order of inclusion and satisfies the following property: for any $\frakp,\frakq\in\opn{Spec}(R)$, $\{\frakp\}\leq \{\frakq\}$ if and only if $\frakp\leq \frakq$ in the fixed linear extension. This is the colexicographic order, which is already implemented in Mathlib under the name \lean{Finset.Colex}. Finally, we take the Dedekind-MacNeille completion of this colexicographically ordered set, which makes it into a complete linear order.
\begin{leancode}[hyperurl={https://github.com/YijunYuan/HarderNarasimhan/blob/c6eb8349df607ee99fa2cab62aee7d9ec367418a/HarderNarasimhan/Coprimary/Defs.lean\#L100}]
DedekindCut (Colex (Finset (LinearExtension (PrimeSpectrum R))))
\end{leancode}

\subsubsection{The pay-off function $\mu_M$}
For any submodules $N_1\subsetneq N_2$ of $M$, the value of the pay-off function $\mu_M$ at $(N_1,N_2)$ is defined by: 
$$\mu_M(N_1,N_2)\coloneqq \opn{Ass}(N_2/N_1).$$

In this project, we first define the set of associated primes of the subquotient $N_2/N_1$ as follows:
\begin{leancode}[hyperurl={https://github.com/YijunYuan/HarderNarasimhan/blob/c6eb8349df607ee99fa2cab62aee7d9ec367418a/HarderNarasimhan/Coprimary/Defs.lean\#L69-L73}]
def subquotientAssociatedPrimes (I : StrictIntvl (Submodule R M)) :
    Set (LinearExtension (PrimeSpectrum R)) :=
  {q | q.asIdeal ∈ associatedPrimes R (I.right ⧸ I.left.submoduleOf I.right)}
\end{leancode}
Noetherianity and finite generation make this set finite; the implementation uses the
Mathlib theorem \lean{associatedPrimes.finite}. After taking \lean{.toFinset}, this set is embedded into the complete order via
\lean{toColex} and \lean{DedekindCut.principal}.  The result is bundled as a
\lean{PayoffFunction}, so all general operations can be used with dot notation:
\begin{leancode}[hyperurl={https://github.com/YijunYuan/HarderNarasimhan/blob/c6eb8349df607ee99fa2cab62aee7d9ec367418a/HarderNarasimhan/Coprimary/Defs.lean\#L95-L101}]
def payoff :
    PayoffFunction (Submodule R M)
      (DedekindCut (Colex (Finset (LinearExtension (PrimeSpectrum R))))) :=
  ⟨fun I ↦ .principal (toColex (subquotientAssociatedPrimes I).toFinset)⟩
\end{leancode}

\subsubsection{Calculation of $\mu_{M,A}$}\label{sec:mua}
In this Harder-Narasimhan game, there is an explicit characterization of $\mu_{M,A}$, which is frequently used in the proof of \Cref{thm:cofil}:
\begin{proposition}[cf. {\cite[Proposition 3.12]{chen2023hardernarasimhangames}}]\label{prop:mua}
    Let $N_1\subsetneq N_2$ be two submodules of $M$, and let $\frakp$ be the least element of $\mu_M(N_1,N_2)$ in the totally ordered set $\opn{Spec}(R)$. Then one has $\mu_{M,A}(N_1,N_2)=\{\frakp\}$.
\end{proposition}
\begin{leancode}[hyperurl={https://github.com/YijunYuan/HarderNarasimhan/blob/c6eb8349df607ee99fa2cab62aee7d9ec367418a/HarderNarasimhan/Coprimary/Semistability.lean\#L293-L321}]
lemma A_payoff (I : StrictIntvl (Submodule R M)) :
    (payoff R M).A I =
      .principal (toColex {(subquotientAssociatedPrimes I).toFinset.min'
        (subquotientAssociatedPrimes_nonempty I)}) := by
\end{leancode}
The optimal first move is constructed by localizing the subquotient at this least prime and
lifting the kernel of the localization map back to the interval. The commutative-algebra input is the following result:
\begin{proposition}[{cf. \cite[Ch IV, §1, n°2, Proposition 6]{bourbakiAlgebreCommutative2006a}}]\label{prop:bourbaki}
Let $S$ be a submonoid of a commutative Noetherian ring $R$.  The associated primes of the
quotient of $M$ by the kernel of $M\to S^{-1}M$ are precisely those associated primes of $M$
which are disjoint from $S$.
\end{proposition}
\begin{leancode}[hyperurl={https://github.com/YijunYuan/HarderNarasimhan/blob/c6eb8349df607ee99fa2cab62aee7d9ec367418a/HarderNarasimhan/Coprimary/AssociatedPrimes.lean\#L119-L133}]
theorem associatedPrimes_quot_ker_mkLinearMap [IsNoetherianRing R] :
    associatedPrimes R (M ⧸ LinearMap.ker (LocalizedModule.mkLinearMap S M)) =
      { p ∈ associatedPrimes R M | p.carrier ∩ S = ∅ }
\end{leancode}
\begin{proof}
  Put $K=\ker(M\to S^{-1}M)$ and $Q=M/K$. Every element of $K$ is annihilated
  by some element of $S$, so every prime in $\opn{Ass}(K)$ meets $S$.
  On the other hand, $Q$ embeds into $S^{-1}M$, on which every element of $S$
  acts invertibly. Multiplication by any element of $S$ is therefore injective
  on $Q$, so every prime in $\opn{Ass}(Q)$ is disjoint from $S$.

  Let $\frakp\in\opn{Ass}(Q)$. Since $\frakp\cap S=\varnothing$, every element
  of $S$ becomes a unit in $R_\frakp$. Thus $K_\frakp=0$, and exactness of
  localization gives $Q_\frakp\cong M_\frakp$. The localization theorem for
  associated primes yields
  $$\frakp R_\frakp\in\opn{Ass}_{R_\frakp}(Q_\frakp)
    =\opn{Ass}_{R_\frakp}(M_\frakp).$$
  Since $R$ is Noetherian, contracting this associated prime to $R$ gives
  $\frakp\in\opn{Ass}(M)$. Together with the disjointness established above,
  this proves one inclusion.

  Conversely, the exact sequence $0\to K\to M\to Q\to0$ gives
  $$\opn{Ass}(M)\subseteq\opn{Ass}(K)\cup\opn{Ass}(Q).$$
  If $\frakp\in\opn{Ass}(M)$ is disjoint from $S$, it cannot belong to
  $\opn{Ass}(K)$ by the first observation. Hence $\frakp\in\opn{Ass}(Q)$,
  proving the reverse inclusion.
\end{proof}

\begin{remark}\label{rmk:62533}
  The proof of \Cref{prop:mua} follows \cite[Proposition 3.12]{chen2023hardernarasimhangames}, with a slight modification to the proof of the lower bound $\{\frakp\}\leq\mu_{M,A}(N_1,N_2)$: for $Q=N_2/N_1$, we work directly with minimal primes over $\opn{Ann}_R(Q)$, which belong to $\opn{Ass}(Q)$, rather than passing through the minimal elements of $\opn{Supp}(Q)$. This simplifies the formalization by using Mathlib's existing results directly and avoiding auxiliary lemmas about support.

  In the initial version of this project, we used \cite[Tag 02CE]{stacks-project}, which identifies the minimal associated primes with the minimal elements of the support, and formalized the full statement of \cite[Ch IV, §1, n°2, Proposition 6]{bourbakiAlgebreCommutative2006a}. Subsequent code updates and refactoring showed that both statements were stronger than needed here. The current proof uses Mathlib's existing results on minimal primes over annihilators without explicitly invoking \cite[Tag 02CE]{stacks-project}, and retains only Bourbaki's quotient formula, stated in \Cref{prop:bourbaki}. See \Cref{sec:53626} for an account of the original formalizations.
\end{remark}
\begin{remark}
    There is an inconsistency in the design of certain commutative-algebra-related classes in Mathlib. Given a ring $R$, among the various subsets of its prime spectrum, some are implemented as a \lean{Set} of type \lean{PrimeSpectrum R}:
\begin{leancode}[hyperurl={https://github.com/leanprover-community/mathlib4/blob/db584cd6d46c92f209a44c0f1c829460d327499d/Mathlib/RingTheory/Support.lean\#L47-L50}]
/-- The support of a module, defined as the set of primes `p` such that `Mₚ ≠ 0`. -/
def Module.support : Set (PrimeSpectrum R) :=
  { p | Nontrivial (LocalizedModule p.asIdeal.primeCompl M) }
\end{leancode}
\noindent while others are implemented as a \lean{Set} of type \lean{Ideal R}:
\begin{leancode}[hyperurl={https://github.com/leanprover-community/mathlib4/blob/db584cd6d46c92f209a44c0f1c829460d327499d/Mathlib/RingTheory/Ideal/AssociatedPrime/Basic.lean\#L96-L104}]
/-- `IsAssociatedPrime I M` if the prime ideal `I` is the radical of the annihilator
of some `x : M`. -/
def IsAssociatedPrime : Prop :=
  (⊥ : Submodule R M).IsAssociatedPrime I

variable (R) in
/-- The set of associated primes of a module. -/
def associatedPrimes : Set (Ideal R) :=
  { I | IsAssociatedPrime I M }
\end{leancode}
This inconsistency causes extra work when one needs to compare these different subsets of $\opn{Spec}(R)$. For example, 
\begin{leancode}
associatedPrimes R M ⊆ Module.support R M
\end{leancode}
\noindent is not valid Lean code for the well-known inclusion $\opn{Ass}(M)\subseteq \opn{Supp}(M)$. Instead, one has to write the following
\begin{leancode}
∀ p, (hp : p ∈ associatedPrimes R M) → ⟨p, hp.1⟩ ∈ Module.support R M 
\end{leancode}
\noindent which is less concise.
\end{remark}

\subsection{Proof of the coprimary-filtration theorem}
\cite[Theorem 3.15]{chen2023hardernarasimhangames} states, without providing details, that \Cref{thm:cofil} ``is a direct consequence'' of \Cref{thm:hnthm}. While this strategy is correct, considerable detail is required to obtain a complete formalization.

We begin with the following characterization of the semistability of the Harder-Narasimhan game $(\scrL_M,S_R,\mu_M)$:
\begin{lemma}[cf. {\cite[Remark 3.14]{chen2023hardernarasimhangames}}]\label{lem:semistable_coprimary}
    The Harder-Narasimhan game $(\scrL_M,S_R,\mu_M)$ is semistable if and only if $M$ is coprimary.
\begin{leancode}[hyperurl={https://github.com/YijunYuan/HarderNarasimhan/blob/c6eb8349df607ee99fa2cab62aee7d9ec367418a/HarderNarasimhan/Coprimary/Semistability.lean\#L364-L412}]
theorem isSemistable_iff_existsUnique_associatedPrime [Nontrivial M] :
    (payoff R M).IsSemistable ↔ ∃! p, p ∈ associatedPrimes R M
\end{leancode}
\end{lemma}
\begin{proof}
  Let $\frakp_0$ be the least element of the nonempty finite set $\opn{Ass}(M)$
  in the fixed total order on $\opn{Spec}(R)$. For every nonzero submodule $N$ of $M$,
  the set $\opn{Ass}(N)$ is also nonempty and finite; denote its least element by
  $\frakp_N$. Since $\opn{Ass}(N)\subseteq\opn{Ass}(M)$, one has
  $\frakp_0\leq\frakp_N$. By \Cref{prop:mua} and the order on singleton sets,
  $$\mu_{M,A}(0,M)=\{\frakp_0\}\leq\{\frakp_N\}=\mu_{M,A}(0,N).$$
  On the other hand, semistability means that
  $\mu_{M,A}(0,N)\not>\mu_{M,A}(0,M)$ for every such $N$.
  Since $S_R$ is totally ordered, the game is therefore semistable if and only if
  $$\mu_{M,A}(0,N)=\mu_{M,A}(0,M)=\{\frakp_0\}
    \qquad\text{for every nonzero submodule }N\subseteq M.$$

  Suppose first that the game is semistable, and let $\frakp\in\opn{Ass}(M)$.
  Choose $x\in M$ whose annihilator is $\frakp$. The cyclic submodule $Rx$ is
  nonzero and isomorphic to $R/\frakp$, so
  $\opn{Ass}(Rx)=\{\frakp\}$. Applying \Cref{prop:mua} and the equality above gives
  $$\{\frakp\}=\mu_{M,A}(0,Rx)=\mu_{M,A}(0,M)=\{\frakp_0\}.$$
  Thus every associated prime of $M$ equals $\frakp_0$, and $M$ is coprimary.

  Conversely, suppose that $\opn{Ass}(M)=\{\frakp_0\}$. For every nonzero
  submodule $N\subseteq M$, the inclusion
  $\varnothing\neq\opn{Ass}(N)\subseteq\opn{Ass}(M)$ forces
  $\opn{Ass}(N)=\{\frakp_0\}$. Hence \Cref{prop:mua} gives
  $\mu_{M,A}(0,N)=\{\frakp_0\}=\mu_{M,A}(0,M)$ for every such $N$,
  which proves semistability.
\end{proof}

The existence of a coprimary filtration of $M$ is proved by showing that the Harder-Narasimhan filtration associated with the game $(\scrL_M,S_R,\mu_M)$ is a coprimary filtration of $M$:
\begin{lemma}\label{lem:hn_coprimary}
    Let $\calF_{\opn{HN}}$ be the (unique) Harder-Narasimhan filtration of length $n$ associated to the game $(\scrL_M,S_R,\mu_M)$. Then
    \begin{enumerate}
        \item For $i=0,\cdots,n-1$, $\calF_{\opn{HN}}(i+1)/\calF_{\opn{HN}}(i)$ is coprimary.
\begin{leancode}[hyperurl={https://github.com/YijunYuan/HarderNarasimhan/blob/c6eb8349df607ee99fa2cab62aee7d9ec367418a/HarderNarasimhan/Coprimary/Filtration.lean\#L49-L60}]
lemma piecewise_isCoprimary
    {M : Type*} [AddCommGroup M] [Module R M] [Module.Finite R M]
    (F : (Coprimary.payoff R M).HarderNarasimhanFiltration) :
    ∀ i < F.length, IsCoprimary R (F (i + 1) ⧸ (F i).submoduleOf (F (i + 1)))
\end{leancode}
        \item Suppose $\calF_{\opn{HN}}(i+1)/\calF_{\opn{HN}}(i)$ is $\frakp_i$-coprimary for $i=0,\cdots,n-1$. Then the sequence $\{\frakp_i\}$ is decreasing in the total order on $\opn{Spec}(R)$.
\begin{leancode}[hyperurl={https://github.com/YijunYuan/HarderNarasimhan/blob/c6eb8349df607ee99fa2cab62aee7d9ec367418a/HarderNarasimhan/Coprimary/Filtration.lean\#L79}]
∀ (i : ℕ),
  i + 1 < FHN.length →
    ∀ (p q : PrimeSpectrum R),
      p.asIdeal ∈ associatedPrimes R (↥(FHN (i + 2)) ⧸ (FHN (i + 1)).submoduleOf (FHN (i + 2))) →
        q.asIdeal ∈ associatedPrimes R (↥(FHN (i + 1)) ⧸ (FHN i).submoduleOf (FHN (i + 1))) →
          toLinearExtension p < toLinearExtension q
\end{leancode}
    \end{enumerate}
\end{lemma}
\begin{proof}
  For every integer $i=0,\cdots,n-1$, let $Q_i\coloneqq \calF_{\opn{HN}}(i+1)/\calF_{\opn{HN}}(i)$, which is a non-zero finitely generated $R$-module. The submodule correspondence gives an order isomorphism
  $$[\calF_{\opn{HN}}(i),\calF_{\opn{HN}}(i+1)]\longrightarrow\scrL_{Q_i},\qquad W\longmapsto W/\calF_{\opn{HN}}(i).$$
  For $\calF_{\opn{HN}}(i)\subseteq U\subsetneq V\subseteq \calF_{\opn{HN}}(i+1)$, the canonical isomorphism
  $$(V/\calF_{\opn{HN}}(i))/(U/\calF_{\opn{HN}}(i))\cong V/U$$
  identifies the associated primes, hence the pay-off functions, of the corresponding
  strict intervals. Consequently, the restricted game on $[\calF_{\opn{HN}}(i),\calF_{\opn{HN}}(i+1)]$ is semistable
  if and only if the game $(\scrL_{Q_i},S_R,\mu_{Q_i})$ is semistable. The former is
  semistable by the definition of a Harder-Narasimhan filtration. And the first assertion follows from
  \Cref{lem:semistable_coprimary}.

  Write $\opn{Ass}(Q_i)=\{\frakp_i\}$. By \Cref{prop:mua},
  $$\mu_{M,A}(\calF_{\opn{HN}}(i),\calF_{\opn{HN}}(i+1))=\{\frakp_i\}.$$
  For $0\leq i<n-1$, the second defining property of a Harder-Narasimhan filtration gives
  $$\mu_{M,A}(\calF_{\opn{HN}}(i),\calF_{\opn{HN}}(i+1))\not\leq\mu_{M,A}(\calF_{\opn{HN}}(i+1),\calF_{\opn{HN}}(i+2)).$$
  Since $S_R$ is totally ordered, this means $\{\frakp_{i+1}\}<\{\frakp_i\}$ and consequently $\frakp_{i+1}<\frakp_i$ by the definition of the colexicographic order. This proves the second assertion.
\end{proof}
Conversely, to deduce the uniqueness of the coprimary filtration from the uniqueness of the Harder-Narasimhan filtration, we need to verify that every coprimary filtration is a Harder-Narasimhan filtration of $(\scrL_M,S_R,\mu_M)$:
\begin{lemma}\label{lem:coprimary_hn}
    Let $\calF_{\opn{CP}}$ be a coprimary filtration of $M$ satisfying the decreasing-prime condition of \Cref{thm:cofil}. Denote by $m$ the length of $\calF_{\opn{CP}}$.  Then
    \begin{enumerate}
        \item For $i=0,\cdots,m-1$, the restricted pay-off function $\mu_{M,[\calF_{\opn{CP}}(i),\calF_{\opn{CP}}(i+1)]}$ is semistable.
\begin{leancode}[hyperurl={https://github.com/YijunYuan/HarderNarasimhan/blob/c6eb8349df607ee99fa2cab62aee7d9ec367418a/HarderNarasimhan/Coprimary/Filtration.lean\#L112}]
∀ (i : ℕ) (hi : i < FCP.length),
  ((Coprimary.payoff R M).restrict { left := FCP i, right := FCP (i + 1), lt := ⋯ }).IsSemistable
\end{leancode}
        \item For $i=1,\cdots,m-1$, one has
        $$\mu_{M,A}(\calF_{\opn{CP}}(i-1),\calF_{\opn{CP}}(i))>\mu_{M,A}(\calF_{\opn{CP}}(i),\calF_{\opn{CP}}(i+1)).$$
\begin{leancode}[hyperurl={https://github.com/YijunYuan/HarderNarasimhan/blob/c6eb8349df607ee99fa2cab62aee7d9ec367418a/HarderNarasimhan/Coprimary/Filtration.lean\#L118}]
∀ (i : ℕ) (hi : i + 1 < FCP.length),
  ¬(Coprimary.payoff R M).A { left := FCP i, right := FCP (i + 1), lt := ⋯ } ≤
      (Coprimary.payoff R M).A { left := FCP (i + 1), right := FCP (i + 2), lt := ⋯ }
\end{leancode}
    \end{enumerate}
\end{lemma}
\begin{proof}
  For every integer $i=0,\cdots,m-1$, let
  $Q_i\coloneqq\calF_{\opn{CP}}(i+1)/\calF_{\opn{CP}}(i)$, which is a non-zero finitely generated coprimary $R$-module.
  Hence \Cref{lem:semistable_coprimary} implies that
  $(\scrL_{Q_i},S_R,\mu_{Q_i})$ is semistable.
  As in the proof of \Cref{lem:hn_coprimary}, the submodule correspondence
  $$[\calF_{\opn{CP}}(i),\calF_{\opn{CP}}(i+1)]\longrightarrow\scrL_{Q_i},
    \qquad W\longmapsto W/\calF_{\opn{CP}}(i)$$
  identifies the pay-off functions on corresponding strict intervals, and therefore
  preserves semistability. Thus $\mu_{M,[\calF_{\opn{CP}}(i),\calF_{\opn{CP}}(i+1)]}$
  is semistable, proving the first assertion.

  Write $\opn{Ass}(Q_i)=\{\frakp_i\}$. By the defining order condition on a
  coprimary filtration, $\frakp_{i-1}>\frakp_i$ for $i=1,\cdots,m-1$.
  By \Cref{prop:mua}, one has
  $$\mu_{M,A}(\calF_{\opn{CP}}(i),\calF_{\opn{CP}}(i+1))=\{\frakp_i\},$$
  and consequently
  $$\mu_{M,A}(\calF_{\opn{CP}}(i-1),\calF_{\opn{CP}}(i))
    =\{\frakp_{i-1}\}>\{\frakp_i\}
    =\mu_{M,A}(\calF_{\opn{CP}}(i),\calF_{\opn{CP}}(i+1)).$$
\end{proof}

\begin{proof}[Completion of the proof of \Cref{thm:cofil}]
  We first check the hypotheses of \Cref{thm:hnthm}. Since $M$ is Noetherian,
  $\scrL_M$ satisfies the ascending chain condition. The second isomorphism theorem
  gives $X/(X\cap Y)\cong(X+Y)/Y$, hence
  $$\mu_M(X\cap Y,X)=\mu_M(Y,X+Y)$$
  whenever $X\not\subseteq Y$; in particular, $\mu_M$ is convex.
  For the $\mu_{M,A}$-descending chain condition, consider a descending chain
  $L_0\supsetneq L_1\supsetneq\cdots\supsetneq N$. By \Cref{prop:mua}, strictly
  increasing values $\mu_{M,A}(N,L_j)$ would give infinitely many distinct least
  primes of $\opn{Ass}(L_j/N)$. This is impossible, since all these primes belong
  to the finite set $\opn{Ass}(L_0/N)$. Finally, $S_R$ is totally ordered, so the
  admissibility hypothesis holds as well. Thus \Cref{thm:hnthm} and
  \Cref{lem:hn_coprimary,lem:coprimary_hn} give the existence and uniqueness of
  the required coprimary filtration.

  Write this filtration as $0=M_0\subsetneq\cdots\subsetneq M_n=M$, with
  $\opn{Ass}(M_{i+1}/M_i)=\{\frakp_i\}$. Applying the associated-prime inclusion
  for the exact sequences $0\to M_i\to M_{i+1}\to M_{i+1}/M_i\to0$ inductively gives
  $$\opn{Ass}(M)\subseteq\bigcup_{i=0}^{n-1}\opn{Ass}(M_{i+1}/M_i)
    =\{\frakp_0,\ldots,\frakp_{n-1}\}.$$
  For the reverse inclusion, the canonical Harder-Narasimhan filtration satisfies
  $$\mu_{M,A}(0,M_{i+1})=\mu_{M,A}(M_i,M_{i+1})=\{\frakp_i\}.$$
  The first equality is
  \lean{PayoffFunction.hnFiltration_A_bot_eq_A}\link{https://github.com/YijunYuan/HarderNarasimhan/blob/c6eb8349df607ee99fa2cab62aee7d9ec367418a/HarderNarasimhan/Filtration/Exists.lean\#L172-L189};
  the second follows from \Cref{prop:mua}. Applying \Cref{prop:mua} also to
  $(0,M_{i+1})$ shows that $\frakp_i$ is the least element of
  $\opn{Ass}(M_{i+1})\subseteq\opn{Ass}(M)$. Hence every $\frakp_i$ belongs to
  $\opn{Ass}(M)$, proving the asserted equality.
\end{proof}

\section{Application 2: the Jordan-Hölder filtration}
For every non-trivial semistable vector bundle $\calE$ over a smooth projective curve, there exists a filtration of $\calE$ by sub-bundles
$$0=\calE_0\subsetneq \calE_1\subsetneq \cdots\subsetneq \calE_n=\calE$$
such that each quotient $\calE_{i+1}/\calE_i$ is stable with the same slope as $\calE$. Such a filtration is called a \textbf{Jordan-Hölder filtration} or \textbf{stable filtration} of $\calE$ (cf. \cite[Section 1.5]{Huybrechts_Lehn_2010}). Although the Jordan-Hölder filtration is generally not unique, the length of such a filtration is always the same.

In the setting of \cite{chen2023hardernarasimhangames}, the Jordan-Hölder filtration is formulated as follows:
\begin{definition}
  A \textbf{Jordan-Hölder filtration} of a Harder-Narasimhan game $(\scrL,S,\mu)$ is a decreasing sequence
  $$\top=y_0>y_1>\cdots>y_n=\bot$$
  such that, for any $i\in\{1,\cdots,n\}$, $\mu(y_i,y_{i-1})=\mu(\bot,\top)$ and for every $z\in \scrL$ such that $y_i<z<y_{i-1}$, one has $\mu(y_i,z)<\mu(y_i,y_{i-1})$.
\end{definition}
\begin{leancode}[hyperurl={https://github.com/YijunYuan/HarderNarasimhan/blob/c6eb8349df607ee99fa2cab62aee7d9ec367418a/HarderNarasimhan/JordanHolder/Defs.lean\#L110-L136}]
structure JordanHolderFiltration (μ : PayoffFunction ℒ S) where
  toFun : ℕ → ℒ
  length : ℕ
  antitone : Antitone toFun
  head_eq_top : toFun 0 = ⊤
  length_eq_bot : toFun length = ⊥
  strictAntiOn : StrictAntiOn toFun (Set.Iic length)
  step_payoff_eq : ∀ i, (hi : i < length) →
    μ ⟨toFun (i + 1), toFun i, strictAntiOn hi.le hi (lt_add_one i)⟩ = μ ⊤
  payoff_lt_of_between : ∀ i, (hi : i < length) → ∀ z : ℒ, (h' : toFun (i + 1) < z) →
    z < toFun i →
    μ ⟨toFun (i + 1), z, h'⟩ <
      μ ⟨toFun (i + 1), toFun i, strictAntiOn hi.le hi (lt_add_one i)⟩
\end{leancode}

The following theorem provides sufficient conditions for the existence of a Jordan-Hölder filtration:
\begin{theorem}[cf. {\cite[Theorem 4.25]{chen2023hardernarasimhangames}}]\label{thm:jhfil}
  Let $(\scrL,S,\mu)$ be a semistable Harder-Narasimhan game. Then it admits a Jordan-Hölder filtration if the following conditions are satisfied:
  \begin{enumerate}
    \item $\scrL$ is a bounded lattice that satisfies the ascending chain condition;
    \item the order on $S$ is total;
    \item the pay-off function $\mu$ is slope-like and satisfies the eventually-$\top$ descending chain condition; 
    \item $\mu(\bot,\top)\neq\top$, the
      greatest element of $S$. This is wrapped into a \lean{Prop}-valued \lean{class} \lean{FiniteTotalPayoff}.
  \end{enumerate}
\end{theorem}

The construction is greedy.  Starting at $\top$, it chooses a maximal $p$ below the current
term for which $\mu(\bot,p)=\mu(\bot,\top)$; if no such $p$ exists, it moves to $\bot$:
\begin{leancode}[hyperurl={https://github.com/YijunYuan/HarderNarasimhan/blob/c6eb8349df607ee99fa2cab62aee7d9ec367418a/HarderNarasimhan/JordanHolder/Exists.lean\#L46-L59}]
def JHFil (k : ℕ) : ℒ :=
  match k with
  | 0 => ⊤
  | n + 1 =>
    let 𝒮 := {p : ℒ | ∃ h : ⊥ < p, p < JHFil n ∧ μ ⟨⊥, p, h⟩ = μ ⊤}
    if h𝒮 : 𝒮.Nonempty then
      (hacc.wf.has_min 𝒮 h𝒮).choose
    else
      ⊥
\end{leancode}
The seesaw property supplies the step pay-off and stability conditions, while
\lean{EventuallyTopDCC} and \lean{FiniteTotalPayoff} force termination.  The result is
packaged into a \lean{Nonempty} instance, and we also provide a \lean{RelSeries} variant of this result:
\begin{leancode}[hyperurl={https://github.com/YijunYuan/HarderNarasimhan/blob/c6eb8349df607ee99fa2cab62aee7d9ec367418a/HarderNarasimhan/JordanHolder/Defs.lean\#L230-L236}]
def jordanHolderRel (μ : PayoffFunction ℒ S) : SetRel ℒ ℒ :=
  {(x, y) | ∃ h : y < x, μ ⟨y, x, h⟩ = μ ⊤ ∧
    ∀ z : ℒ, (h' : y < z) → z < x → μ ⟨y, z, h'⟩ < μ ⟨y, x, h⟩}
\end{leancode}
\begin{leancode}[hyperurl={https://github.com/YijunYuan/HarderNarasimhan/blob/c6eb8349df607ee99fa2cab62aee7d9ec367418a/HarderNarasimhan/JordanHolder/Exists.lean\#L195-L203}]
theorem exists_relSeries_jordanHolderRel :
    ∃ s : RelSeries (μ.jordanHolderRel), s.head = ⊤ ∧ s.last = ⊥ 
\end{leancode}

\begin{remark}
In contrast to \lean{HarderNarasimhanFiltration}, the implementation deliberately provides no
\lean{Inhabited} instance: a Jordan-Hölder filtration need not be unique, and the greedy
construction is not exposed as a canonical public choice.
\end{remark}

\begin{remark}
The condition ``for any $i\in\{1,\cdots,n\}$ and $z\in \scrL$ such that $y_i<z<y_{i-1}$,
one has $\mu(y_i,z)<\mu(y_i,y_{i-1})$'' is the field
\lean{payoff_lt_of_between}.  The following theorem shows that it is equivalent to the stability
of the restriction of $\mu$ to every step interval:
\begin{leancode}[hyperurl={https://github.com/YijunYuan/HarderNarasimhan/blob/c6eb8349df607ee99fa2cab62aee7d9ec367418a/HarderNarasimhan/JordanHolder/Stability.lean\#L148-L155}]
theorem piecewise_isStable_iff (hsa : ∀ i j : ℕ, i < j → j ≤ n → f j < f i) :
    (∀ i : ℕ, (hi : i < n) →
        (μ.restrict ⟨f (i + 1), f i, hsa i (i + 1) (lt_add_one i) hi⟩).IsStable) ↔
      ∀ i : ℕ, (hi : i < n) → ∀ z : ℒ, (h' : f (i + 1) < z) → z < f i →
        μ ⟨f (i + 1), z, h'⟩ < μ ⟨f (i + 1), f i, hsa i (i + 1) (lt_add_one i) hi⟩
\end{leancode}
Therefore, the restriction of a semistable Harder-Narasimhan game to each interval between consecutive terms of a Jordan-Hölder filtration is stable.
\end{remark}
\begin{proof}
  Let $(\scrL,S,\mu)$ be as in \Cref{thm:jhfil}, and let
  $f\colon\bfN\to\scrL$ satisfy $f(0)>\cdots>f(n)$.
  Fix $i<n$, put $a=f(i+1)$, $b=f(i)$, and $c=\mu(a,b)$, and work with
  the restricted Harder-Narasimhan game $(\scrL_{[a,b]},S,\mu_{[a,b]})$.
  The chain conditions and slope-like property pass to restrictions, so
  $\mu_A=\mu_{\min}$ on subintervals,
  and the Nash-equilibrium characterization identifies semistability on $[a,b]$
  with $\mu_{\min}(a,b)=\mu_{\max}(a,b)=c$.

  If $\mu(a,z)<c$ for every $a<z<b$, the seesaw property gives
  $\mu(z,b)>c$. Hence the restricted game is semistable, and
  $$\mu_A(a,z)=\mu_{\min}(a,z)\leq\mu(a,z)<c=\mu_A(a,b)
    \qquad(a<z<b),$$
  proving stability.

  Conversely, stability gives $\mu(a,z)\leq c$ and
  $\mu_{\min}(a,z)=\mu_A(a,z)<c$ for $a<z<b$.
  If $\mu(a,z)=c$, semistability gives $\mu(a,u)\leq c$ for $a<u<z$,
  so the seesaw property implies $\mu(u,z)\geq c$.
  Including $u=a$ yields $\mu_{\min}(a,z)=c$, a contradiction.
  Thus $\mu(a,z)<c$, as required.
\end{proof}

As with the classical Jordan-Hölder filtration, under certain conditions the length of a Jordan-Hölder filtration is unique:
\begin{theorem}[cf. {\cite[Remark 4.26]{chen2023hardernarasimhangames}}]\label{thm:jhlen}
  Let $(\scrL,S,\mu)$ be a semistable Harder-Narasimhan game as in \Cref{thm:jhfil}. Suppose further that
  \begin{enumerate}
  \item The lattice $\scrL$ is \textbf{modular}, i.e. for any $x,y,z\in\scrL$ such that $x\leq z$, one has
  $$x\vee (y\wedge z) = (x\vee y)\wedge z.$$
  \item The pay-off function $\mu$ is affine.
\end{enumerate}
  Then any two Jordan-Hölder filtrations of $(\scrL,S,\mu)$ have the same length.
\end{theorem}
\begin{leancode}[hyperurl={https://github.com/YijunYuan/HarderNarasimhan/blob/c6eb8349df607ee99fa2cab62aee7d9ec367418a/HarderNarasimhan/JordanHolder/Length.lean\#L576-L582}]
theorem JordanHolderFiltration.length_eq (F G : μ.JordanHolderFiltration) :
    F.length = G.length
\end{leancode}
Here the ambient instances include \lean{IsModularLattice ℒ} and \lean{μ.IsAffine}, in
addition to the existence hypotheses above.  

The proof of \Cref{thm:jhlen} follows the inductive framework of
\cite[Remark 4.26]{chen2023hardernarasimhangames}:
\begin{quotation}
\itshape Let $\{x_i\}_{i\geq 0}$ (resp. $\{y_j\}_{j\geq 0}$) be a Jordan-Hölder filtration of the game $(\scrL,S,\mu)$ of length $m$ (resp. $n$). We prove by construction that the restriction of this game to the interval $\scrL_{[x_{m-1},\top]}$ admits a Jordan-Hölder filtration of length not greater than $n-1$. Iterating this procedure we obtain $m\leq n$. Then by symmetry we obtain $m=n$.
\end{quotation}
The iterative procedure described above is a typical situation that calls for induction. However, one cannot directly ``apply induction on the length of the Jordan-Hölder filtration of a fixed game'', as is commonly done in informal mathematical writing, because the induction hypothesis concerns the Jordan-Hölder filtration of a different game, namely the restriction of the original game to a proper subinterval. The appropriate way to structure the induction is therefore to formulate the following statement:
\begin{lemma}
  For any $n\in \bfN$ and any semistable Harder-Narasimhan game $(\scrL,S,\mu)$ that satisfies the conditions in \Cref{thm:jhlen}, if it admits a Jordan-Hölder filtration of length not greater than $n$, then every Jordan-Hölder filtration of this game has length not greater than $n$.
\end{lemma}
\begin{proof}[Sketch of proof]
  We prove this statement by induction on $n$. 
  
  Since the length of a Jordan-Hölder filtration is always a positive natural number, the base case $n=0$ is trivial. 
  
  For the induction step, let $n$ be a natural number, and suppose the statement holds for $n$. Let $(\scrL,S,\mu)$ be a semistable Harder-Narasimhan game that satisfies the conditions in \Cref{thm:jhlen}, and suppose it admits a Jordan-Hölder filtration $\{y_j\}_{j\geq 0}$ of length not greater than $n+1$. Let $\{x_i\}_{i\geq 0}$ be a Jordan-Hölder filtration of this game of length $m$. We need to show that $m\leq n+1$. We may assume without loss of generality that $m> 1$, for otherwise we are done.

  From $\{y_j\}_{j\geq 0}$, we construct a Jordan-Hölder filtration of $\mu_{[x_{m-1},\top]}$ with length not greater than $n$, and one checks that the restricted game again satisfies the conditions in \Cref{thm:jhlen}. Since the restriction of $\{x_i\}_{i \geq 0}$ to $\scrL_{[x_{m-1},\top]}$, which has length $m-1$, is also a Jordan-Hölder filtration of $\mu_{[x_{m-1},\top]}$, the induction hypothesis implies that $m-1\leq n$, i.e. $m\leq n+1$.
\end{proof}
The following is a sketch of the Lean code of this lemma:
\begin{leancode}[hyperurl={https://github.com/YijunYuan/HarderNarasimhan/blob/c6eb8349df607ee99fa2cab62aee7d9ec367418a/HarderNarasimhan/JordanHolder/Length.lean\#L528-L566}]
lemma length_le_of_exists_length_le (n : ℕ) :
    ∀ {ℒ : Type*} [Nontrivial ℒ] [Lattice ℒ] [BoundedOrder ℒ]
      [WellFoundedGT ℒ] [IsModularLattice ℒ]
      {S : Type*} [CompleteLinearOrder S] {μ : PayoffFunction ℒ S}
      [μ.FiniteTotalPayoff] [μ.IsSlopeLike] [μ.IsSemistable]
      [μ.EventuallyTopDCC] [μ.IsAffine],
      (∃ F : μ.JordanHolderFiltration, F.length ≤ n) →
      ∀ G : μ.JordanHolderFiltration, G.length ≤ n := by
  induction n with
  | zero =>
    intro ℒ _ _ _ _ _ S _ μ _ _ _ _ _ ⟨F, hF⟩ _
    exact (F.length_pos.not_ge hF).elim
  | succ n ih =>
    intro ℒ _ _ _ _ _ S _ μ hfinite hsl _ _ _ ⟨F, hF⟩ G
    ... -- Omitted for brevity.
\end{leancode}
Once this lemma is established, \Cref{thm:jhlen} follows immediately by symmetry.

\appendix
\section{Influence of artificial intelligence}
Although the main body of the formalization was carried out manually by the author, we benefited from artificial intelligence in the following ways:

\subsection{Accelerating the search for relevant lemmas and definitions}
The Mathlib-specific search engine \href{https://leansearch.net/}{\texttt{LeanSearch}} (cf. \cite{gao-etal-2024-semantic-search}), developed by the AI4Math team at Peking University and driven by a fine-tuned large language model, greatly accelerates the search for relevant definitions and lemmas in Mathlib by supporting natural-language queries.

\subsection{Generating docstrings for the formalized code}
As already mentioned, the project was first completed with Lean v4.27.0 and then upgraded to Lean v4.33.0 and completely refactored. During this period, the capability of large language models improved significantly. In the following, we give a brief account of how AI was used to generate docstrings at the different stages of the project.

\subsubsection*{The project with Lean v4.27.0}
Most of the docstrings at this stage were generated by \texttt{Copilot}\link{https://github.com/features/copilot}, an AI-powered pair programmer developed by Microsoft, in \texttt{Visual Studio Code}, with \texttt{GPT-5.3-Codex} (cf. \cite{GPT53Codex}) as the underlying model. The generation was carried out in Agent mode, with a human-written prompt\link{https://github.com/YijunYuan/HarderNarasimhan/blob/paper/prompt-doc.md} provided.
\begin{remark}
  The file \texttt{HNGame.pdf} mentioned in the prompt is a copy of the preprint \cite{chen2023hardernarasimhangames} on which this project is based. The file was fed to the AI to ensure that the generated comments are consistent with the mathematical content of the paper, and to help the AI understand the mathematical context of the code it was commenting on.
\end{remark}
In prolonged interactions with the AI, we observed a decline in the quality of the generated comments: the AI tended to forget some of the requirements specified in the initial prompt. After additional rounds of refinement and manual adjustment, the quality of the documentation reached a merely satisfactory level.

\subsubsection*{The project with Lean v4.33.0}
We used \texttt{Claude Fable 5} (cf. \cite{Fable5}) and \texttt{GPT-6 Astra} (cf. \cite{GPT6Astra}) extensively to generate docstrings for the refactored code. The official Mathlib docstring style guide\link{https://leanprover-community.github.io/contribute/doc.html} and the electronic version of \cite{chen2023hardernarasimhangames} were fed to the AI. The only constraint we imposed was that the generated docstrings should not contain the numbering of theorems, lemmas, definitions, etc.\ in the original paper, since the numbering might change after the paper is revised. The quality of the generated docstrings improved significantly compared with the previous stage.

The author made an effort to maintain the quality of the documentation by reviewing and refining the AI-generated docstrings for mathematical accuracy and consistency with the Lean declarations and the source paper.

\subsection{Proof generation in early stages}\label{sec:53626}
In the earlier version of the project, two required commutative-algebra results (\cite[Ch IV, §1, n°2, Proposition 6]{bourbakiAlgebreCommutative2006a} and \cite[Tag 02CE]{stacks-project}) were taken from outside
\cite{chen2023hardernarasimhangames}. The following account records how their proofs were generated by AI.
\begin{enumerate}
  \item For \cite[Tag 02CE]{stacks-project}, we instructed \texttt{Copilot} in Agent mode to generate a proof and to verify it by invoking \texttt{lake build}.
  
  After several rounds of an unsupervised generate-verify-revise cycle, Copilot successfully completed the formal proof. Throughout the process, we did not provide any mathematical or programming-related hints whatsoever.
  \item For \cite[Ch IV, §1, n°2, Proposition 6]{bourbakiAlgebreCommutative2006a}, repeating the above workflow did not work well: Copilot gave up, reporting that the proof was beyond the available Mathlib infrastructure.
  
    In a second attempt, we switched Copilot to Plan mode, which prompted \texttt{GPT-5.3-Codex} to first produce a high-level proof outline and then progressively refine it into a complete formal proof. Moreover, we supplied a natural-language proof of the proposition\footnote{To obtain this natural-language proof, we provided \texttt{Gemini 3 Pro} (cf. \cite{Gemini3Pro}) with an electronic version of \cite{bourbakiAlgebreCommutative2006a} and asked it to extract the proof together with any necessary intermediate lemmas as a self-contained narrative. The resulting text was then incorporated into the Copilot prompt for proof generation.
    } and indicated that the formalization should be decomposed into several auxiliary lemmas. With this additional guidance, Copilot successfully produced a complete formal proof after several rounds of a generate--verify--revise cycle.
\end{enumerate}

We remark that, with the continued expansion of Mathlib and the rapid progress of AI models, these early formal proofs are no longer optimal: in the current version of the project, part of the code has been removed as it is now superseded by Mathlib, and the rest has been compressed into much shorter proofs.

\subsection{Code refactoring}\label{sec:refactoring}
After this project was released (cf. \url{https://github.com/YijunYuan/HarderNarasimhan/tree/a2d585b56dfc53a9cadef84ed3595289b4282738}), we realized that the code was not in its best form, for several reasons.
\begin{enumerate}
  \item The organization followed the paper's exposition, with separate statement and implementation layers that introduced redundant wrappers.
  \item Unbundled pay-off functions and repeated hypotheses made declarations verbose and operations on restricted games cumbersome.
  \item Large compound statements, such as the seesaw lemma, were inconvenient to reuse, while some custom order-theoretic constructions could now be replaced by Mathlib.
\end{enumerate}
For these reasons, we decided to refactor the code. After some human-AI\footnote{Multiple models were used in this process: \texttt{GPT-5.4} (cf. \cite{GPT54}), \texttt{GPT-5.5} (cf. \cite{GPT55}), \texttt{GPT-5.6 Sol} (cf. \cite{GPT56}), \texttt{Claude Opus 4.7} (cf. \cite{Opus47}), \texttt{Claude Opus 4.8} (cf. \cite{Opus48}), and \texttt{Claude Fable 5}.} collaborative simplification, we asked \texttt{Claude Fable 5} in \texttt{Claude Code} to prepare a refactoring plan. We provided \texttt{Claude Fable 5} with a human-written prompt\textsuperscript{\href{https://github.com/YijunYuan/HarderNarasimhan/blob/c6eb8349df607ee99fa2cab62aee7d9ec367418a/prompts/refactor.md}{\faIcon{external-link-alt}}} with the following requirements specified:
\begin{enumerate}
  \item Combine statements and proofs in the same files, and design the public interface around useful results, keeping purely auxiliary lemmas internal.
  \item Adopt a modular organization modeled after Mathlib, with four main folders: \texttt{PayoffFunction}, \texttt{Filtration}, \texttt{JordanHolder}, and \texttt{Coprimary}.
  \item Reorganize namespaces and use section variables consistently.
  \item Bundle the pay-off function as \texttt{PayoffFunction} and use dot notation for its associated constructions and properties.
  \item Rewrite docstrings and follow Mathlib's conventions for naming\textsuperscript{\href{https://leanprover-community.github.io/contribute/naming.html}{\faIcon{external-link-alt}}}, formatting\textsuperscript{\href{https://leanprover-community.github.io/contribute/style.html}{\faIcon{external-link-alt}}}, and documentation\textsuperscript{\href{https://leanprover-community.github.io/contribute/doc.html}{\faIcon{external-link-alt}}}.
  \item Review the entire project and record a comprehensive plan in a Markdown file, to be discussed and revised with the author before any code changes.
\end{enumerate}

\texttt{Claude Fable 5} produced a Markdown file\textsuperscript{\href{https://github.com/YijunYuan/HarderNarasimhan/blob/c6eb8349df607ee99fa2cab62aee7d9ec367418a/prompts/REFACTOR_PLAN.md}{\faIcon{external-link-alt}}} containing a detailed refactoring plan of about 700 lines. The author asked the AI to discuss the plan with him item by item, and to revise it in light of his feedback.

Afterwards, \texttt{Claude Fable 5} refactored the code according to the plan, and a subsequent polishing pass was carried out manually by the author. The final result is the current version of the project.

We emphasize that deep human involvement was indispensable throughout this process: proposing the initial refactoring roadmap, discussing the plan with the AI and revising it accordingly, and polishing the code after the refactoring. In an era when AI technology is advancing by leaps and bounds, the sheer volume of formalized code is no longer the bottleneck. Rather, what human mathematicians and formalization researchers should now focus on is a better understanding of the mathematics and active participation in shaping the logical architecture of the code.

\printbibliography
\end{document}